\documentclass{amsart}
\usepackage{graphicx} 
\usepackage{float}
\usepackage{amssymb}		
\usepackage{graphicx}		
\usepackage{hyperref}		

\usepackage{amsmath}
\usepackage{amsthm}
\usepackage{bbm}
\usepackage{enumerate}
\usepackage{tikz}
\usepackage{tikz-cd}
\usetikzlibrary{hobby}

\theoremstyle{plain}
\newtheorem{theorem}{Theorem}[section]
\newtheorem{Cor}[theorem]{Corollary}
\newtheorem{lemma}[theorem]{Lemma}
\newtheorem{prop}[theorem]{Proposition}
\newtheorem{Def}[theorem]{Definition}

\numberwithin{theorem}{subsection}
\newtheorem{Ex}[theorem]{Example}

\newtheorem{rthm}{Theorem}

\newenvironment{reptheorem}[1]
  {\renewcommand\therthm{\ref*{#1}}\rthm}
  {\endrthm}
  
\theoremstyle{definition}

\newcommand{\C}{\mathbb{C}}
\newcommand{\R}{\mathbb{R}}
\newcommand{\Q}{\mathbb{Q}}
\newcommand{\Z}{\mathbb{Z}}
\newcommand{\N}{\mathbb{N}}

\newcommand{\Hom}{\text{Hom}}

\newcommand{\unnumberedfootnote}[1]{%
  \begingroup
    \renewcommand{\thefootnote}{}%
    \footnotetext{#1}%
    \addtocounter{footnote}{-1}%
  \endgroup
}
\usepackage{fullpage}				
\usepackage[utf8]{inputenc}			
\usepackage{amsmath, amssymb,amsthm}
\usepackage{graphicx}				
\usepackage{mathtools}
\usepackage{amsfonts}
\usepackage[colorinlistoftodos]{todonotes}	
\usepackage{hyperref}				
\usepackage{tikz}
\usepackage{mathrsfs}
\usepackage{fncylab}
\usepackage{tikz-cd}
\usepackage{wasysym}
\theoremstyle{definition}

\newcommand{\ee}{\epsilon}
\newcommand{\ds}{\displaystyle}
\newcommand{\alp}{\alpha}

\newcommand{\sub}{\subseteq}

\newcommand{\f}{\frac}

\newcommand{\ve}{\varepsilon}

\newcommand{\de}{\delta}

\newcommand{\A}{\mathscr{A}}

\newcommand{\K}{\mathbb{K}}

\newcommand{\vp}{\varphi}
\newcommand{\es}{\emptyset}

\newcommand{\gm}{\gamma}

\newcommand{\br}[1]{\mkern 1.5mu\overline{\mkern-1.5mu#1\mkern-1.5mu}\mkern 1.5mu}

\title{Finitely Summable $K$-homology, the Index Pairing, and Cantor Minimal Systems}
\author{Levi Lorenzo}
\address{Levi Lorenzo, Department of Mathematics, University of Colorado Boulder Campus Box 395, Boulder, CO 80309-0395, USA}
\email{levi.lorenzo@colorado.edu}

\begin{document}

\maketitle
\begin{abstract} We study index pairings for crossed-product $C^*$-algebras arising from minimal actions on the Cantor set. We utilize Putnam's orbit-breaking AF-subalgebras and embeddings to show we can compute any index pairing for Cantor minimal system crossed products using Connes' trace formulas. In the  case of odometers, we show that the associated algebras have uniformly finitely summable $K$-homology.
\end{abstract}

\unnumberedfootnote{2020 \textit{Mathematics Subject Classification} 19K33, 19K56.}
\unnumberedfootnote{This work partially supported by NSF Grants 2000057 and 2247424, PI: Robin Deeley.}
    
\section{Introduction}\label{Chapter:introduction}

 We investigate index pairings between the $K$-theory and $K$-homology of $C^*$-algebras arising from \textit{Cantor minimal systems}. In the analytic model, $K$-homology classes for a $C^*$-algebra are represented by \textit{Fredholm modules}. Fredholm modules are required to satisfy certain conditions modulo the compact operators on a Hilbert space. However, when there is a dense $*$-subalgebra of the $C^*$-algebra on which these relations hold modulo a Schatten $p$-class, we say the Fredholm module is \textbf{$p$-summable} on said subalgebra. If a Fredholm module is $p$-summable for some $p<\infty,$ we say the Fredholm module is \textbf{finitely summable}. When a $K$-homology class has a finitely summable representative, we can use Connes' trace formulas to compute index pairings between this class and $K$-theory classes generated by elements of (matrices over) the specified subalgebra \cite{C}. Further, Connes shows that when a Fredholm module is finitely summable over a dense $*$-subalgebra, we can always enlarge the algebra to one whose matrix amplifications contain representatives of every $K$-theory class \cite{CNDG}.

Thus, if a $C^*$-algebra admits a dense subalgebra on which every $K$-homology class can be represented by a Fredholm module that is $p$-summable for some $p,$ we can compute the pairing between any $K$-homology class and any $K$-theory class over this $C^*$-algebra using a particular trace formula. In this case, we say that the $C^*$-algebra admits \textit{uniformly finitely summable $K$-homology}. Uniform finite summability has been studied for various classes of $C^*$-algebras by Connes \cite{C}, Emerson and Nica \cite{EN}, Goffeng and Mesland \cite{GM}, Crisp \cite{Cr},  Rave \cite{R}, Gerontogiannis \cite{G}, and Puschnigg \cite{Pusch}.

We study summability and the ability to compute index pairings using Connes' trace formulas for $C^*$-algebras arising from \textbf{Cantor minimal systems}. Cantor minimal systems are dynamical systems where $X$ is a Cantor set and $\vp: X \to X$ is a minimal homeomorphism. We form a $C^*$-algebra associated to this system using the crossed-product construction.

We seek to compute index pairings using Connes' trace formulas for such algebras. We use novel methods to show that, when $(X,\vp)$ is a Cantor minimal system, we can compute any index pairing for $C(X) \rtimes_\vp \Z$ using trace formulas. The main results of this paper are: 

 \begin{reptheorem}{eventhm}
    Any index pairing between a class in $K_0(C(X) \rtimes_\vp \Z)$ and one in $K^0(C(X) \rtimes_\vp \Z)$ can be computed using Connes' trace formula for even cycles with $p>0.$
\end{reptheorem}

\begin{reptheorem}{oddthm}
    Let $x \in K^1(C(X) \rtimes_\vp \Z).$ Pairings of $x$ with elements of $\langle[u] \rangle \cong K_1(C(X) \rtimes_\vp \Z)$ can be computed using Connes' trace formula for odd $p>0$ summable cycles.
\end{reptheorem}

The even result is achieved using the orbit-breaking AF-algebras and sequences of Putnam \cite{P}.  Given a Cantor minimal system, $(X,\vp),$ Putnam constructs an AF-subalgebra, $A_{\{y\}},$ of $C(X) \rtimes_\vp \Z$ obtained by ``breaking orbits'' at a point $y.$ Further, Putnam constructs an embedding $\iota:C(X) \rtimes_\vp \Z \to A_{\{y\}},$ see Section \ref{Embed} or \cite[Chapter 6]{P}. These embeddings induce isomorphisms on $K_0$ and, in either order, compose to the identity on $K_0.$ The UCT allows us to prove the analogous results in $K$-homology and show we can compute index pairings using these embeddings.

In the odd case, we show that, for each element in $\text{Hom}(K_1(C(X) \rtimes_\vp \Z), \Z) \cong \Z,$ we can produce a finitely summable odd cycle on $C(X) \rtimes_\vp \Z$ whose image in $\text{Hom}(K_1(C(X) \rtimes_\vp \Z), \Z)$ in the UCT exact sequence is this element. Thus, we have  produced finitely summable cycles to compute any index map with. 

In the case of odometers, we obtain a stronger result. Odometers are precisely the minimal homeomorphisms on the Cantor set that are (metrically) equicontinuous, so we obtain: \begin{reptheorem}{spectripodometer}
Let $\vp: X \to X$ be a metrically equicontinuous action on the Cantor Set (i.e. an odometer, Proposition \ref{equicontinuosodometers}) and also use $\vp:C(X) \to C(X)$ to denote the induced automorphism $f \mapsto f \circ \vp^{-1}.$  Let $\A = \text{span}\{\chi_{C_\mu}| \mu \in Y \}.$ Then, for every $x \in K^*(C(X) \rtimes_\vp \Z)$ and each $p >1,$ there is an unbounded Fredholm module representing $x$ that is $p$-summable on $C_c(\Z,\A).$
\end{reptheorem} 
 We obtain this result by taking the bounded cycles of \cite[Theorem 4.2.1]{R} and lifting them to finitely summable unbounded cycles that exhaust $K^0(C(X)).$ We then use the work of \cite{HSWZ} to take the finitely summable unbounded cycles we constructed to exhaust $K^0(C(X))$ and extend them to finitely summable unbounded cycles that exhaust $K^1(C(X) \rtimes_\vp \Z).$ 

 This work suggests alternate approaches to computing index pairings for a $C^*$-algebra using Connes' trace formulas when we have summable cycles related to the algebra but not necessarily a dense subalgebra on which the summability is uniform.

\section*{Acknowledgments} Many thanks to Robin Deeley for supervising me through my PhD and providing lots of helpful feedback, advice, and support. Thank you to Magnus Goffeng for lots of helpful comments and feedback and to the rest of my PhD committee.
\section{Preliminaries}
    \subsection{Fredholm Modules and $K$-homology}
 The primary reference for analytic $K$-homology is \cite{HR}. $K$-homology classes are represented by \textit{Fredholm modules}. Specifically, $K^0$ is generated by equivalence classes of \textit{graded} Fredholm modules while $K^1$ is generated by equivalence classes of \textit{ungraded} Fredholm modules.
\begin{Def}
Let $A$ be a $C^*$-algebra. An ungraded \textbf{Fredholm module} over $A$ is a triple $(H, \rho,F)$ where: 
\begin{itemize}
    \item $H$ is a separable Hilbert space
    \item $\rho: A \to B(H)$ is a representation
    \item $F \in B(H)$ satisfies:
    \begin{enumerate}
        \item $(F^2-1)\rho(a) \in K(H)$ for each $a \in A,$
        \item $(F-F^*)\rho(a) \in K(H)$ for each $a \in A,$ and
        \item $[F,\rho(a)] \in K(H)$ for each $a \in A.$
    \end{enumerate}
\end{itemize}
\end{Def}
 A \textit{graded Fredholm module} is of the form: \[ (H,\rho,F) = \left(H_+ \oplus H_-, \rho_+ \oplus \rho_-, \begin{bmatrix} 0 & V \\ U & 0 \end{bmatrix} \right)\] where $U-V^*, V-U^*, UV-1, VU-1 \in \K$ and $\rho_+(a) U - U \rho_-(a), V \rho_-(a) - \rho_+(a) V \in K(H)$ for each $a \in A.$ 

We work with crossed product $C^*$-algebras in this thesis, and the $K$-homology of a crossed-product algebra is obtained from the $K$-homology of the base algebra via the Pimsner-Voiculescu exact sequence:
\begin{prop}\label{PVKhom}
For a $C^*$-algebra $A$ and an action $\alp \in \text{Aut}(A)$ there is a Pimsner-Voiculescu exact sequence in K-homology for the crossed product $A \rtimes_\vp \Z:$
\[\begin{tikzcd}
{K^0(A \rtimes_\alp \Z )} \arrow[d] & K^0(A) \arrow[l] & K^0(A) \arrow[l, "1-\alp^*"]   \\
K^1(A) \arrow[r, "1- \alp^*"]   & K^1(A) \arrow[r]  & {K^1(A \rtimes_\alp \Z)} \arrow[u]
\end{tikzcd}\]
\end{prop}
\subsection{Spectral Triples and Unbounded Fredholm modules}
Some of the Fredholm modules in this paper arise from \textit{spectral triples.}
\begin{Def}
An (odd) \textbf{spectral triple} is a triple $(A, H, D)$ where:
\begin{enumerate}[(1)]
    \item $H$ is a separable Hilbert space,
    \item $A$ is a $*$-subalgebra of $B(H),$ and
    \item $D$ is an unbounded self-adjoint operator in $H$ which satisfies, for each $a \in A$:
    \begin{enumerate}
    \item $a \text{Dom}(D) \sub \text{Dom}(D),$
        \item  $[D,a] \in B(H),$ and
        \item $a(1+D^2)^{-1} \in K(H).$
    \end{enumerate}
\end{enumerate}
\end{Def}
If $H$ is equipped with a $\Z/2\Z$ grading, $A$ is an algebra of even operators, and $D$ is odd, then $(A,H,D)$ is an even spectral triple. 

In the case of a $C^*$-algebra, spectral triples determine \textit{unbounded Fredholm modules}:
\begin{Def}
    Suppose $A$ is a $C^*$-algebra, $\A \sub A$ is dense, and $(\A,H,D)$ is a spectral triple. Denote by $\pi$ the action of $\A$ on $H,$ $\pi: \A \to B(H).$ Then, we denote by $(H,\pi, D)$ the \textbf{spectral triple} or \textbf{unbounded Fredholm module} on $\A$ given by $(\A,H,D).$
\end{Def}
There is a notion of summability for spectral triples:
\begin{Def}
    If $A$ is unital and $\text{Tr}(1+D^2)^{-\f p2} < \infty,$ then $(A,H,D)$ is $p-$summable. If such a $p < \infty$ exists, we say $(A,H,D)$ is finitely summable. The infimum of such $p$ is called the \textbf{spectral dimension} of $(A,H,D).$
\end{Def}
Fredholm modules are obtained from spectral triples via the bounded transform: 
\begin{prop}\label{BT}
Let $A$ be a $C^*$-algebra. Suppose $A$ is represented on a Hilbert Space $H$ via $\pi:A \to B(H)$. Further, suppose $\mathscr{A} \sub A$ is dense and $(\pi(\mathscr{A}), H, D)$ is an (even/odd) spectral triple. Then $\left(H, \pi, \frac{D}{(1+D^2)^{\f 12}}\right)$ is an (even/odd) Fredholm module over $A,$ called the \textbf{bounded transform} of $(\pi(\mathscr{A}), H, D).$ In the case that $D$ is invertible, we can use $D|D|^{-1}$ instead of 
$ \frac{D}{(1+D^2)^{\f 12}}.$ 
\end{prop}

\subsection{Index Pairings and Summability}
Our work concerns the index pairing and circumstances in which we can make it more computable. For a more detailed exposition on index pairings and their formulas see \cite{HR}, \cite{C}.
\begin{Def}
    Let $A$ be a separable $C^*$-algebra. Then there is a bilinear map $\langle -, -\rangle_A: K_*(A) \times K^*(A) \to \Z$ called the \textbf{index pairing}.
\end{Def}
The condition on Fredholm modules which facilitates computation of index pairings using trace formulas is \textit{summability}.
\begin{Def}
An odd (even) Fredholm module $(H, \rho, F)$ over a 
$C^*$-algebra $A$ is \textbf{$p$-summable} if there exists a dense sub-algebra $\mathscr{A} \sub A$ such that, for all $a \in \mathscr{A}:$ \[ \rho(a)(F^2-1),\: \rho(a) (F-F^*) \in L^{\f p2}(H), \text{ and } [F,\rho(a)] \in L^p(H).\]
\end{Def}
 Specifically, index pairings with $p$-summable Fredholm modules can be computed using Connes' trace formulas \cite{C}:
\begin{prop}
    Suppose $(H,\rho,F)$ is an odd Fredholm module that is $p$-summable over $\mathscr{A} \sub A$ and $u \in \mathscr{A}$ is a unitary. Then, \[\langle[u], [H, \rho, F] \rangle = \f{(-1)^{\frac{n-1}{2}-1}}{2^{n}}\text{Tr}(\rho(u^*)([F, \rho(u)][F, \rho(u^*)])^{\frac{n-1}{2}} [F, \rho(u)])\] for $n > p$ odd.
\end{prop}
 \begin{prop}
     Suppose $\left (H_+ \oplus H_-,\rho_+ \oplus \rho_- , \begin{bmatrix} 0 & U \\ V & 0 \end{bmatrix}\right)$ is an even Fredholm module that is $p$-summable over $\mathscr{A} \sub A$ and $p$ is a projection in $\mathscr{A}.$ Then,  \begin{footnotesize} \[ \left \langle[p], \left [H_+ \oplus H_-,\rho_+ \oplus \rho_- , \begin{bmatrix} 0 & U \\ V & 0 \end{bmatrix}\right] \right \rangle =  (-1)^{\f{n(n-1)}{2}} \text{Tr}\left(\begin{bmatrix} 1 &0 \\0 & -1 \end{bmatrix} \rho_+ \oplus \rho_-(p)\left(\left[\begin{bmatrix} 0 & U\\ V & 0\end{bmatrix}, \rho_+ \oplus \rho_-(p)\right]\right)^{n}\right) \]\end{footnotesize} for $n > p$ even.
 \end{prop}
In either case, we can extend these formulas to matrix amplifications over $\mathscr{A}.$ Because they are given in terms of a trace, these formulas are often computable. The desire to compute index pairings using trace formulas motivates the study of \textit{uniformly finitely summable} $K$-homology: 
\begin{Def}
A $C^*$- algebra $A$ is said to have \textbf{uniformly $p$-summable} $K$-homology if for $p > 0$ there is a dense $*$-subalgebra $\mathscr{A} \sub A$ such that each class in $K^*(A)$ (for $*=0$ and $1$) can be represented by a Fredholm module which is $p$-summable on $\mathscr{A}.$
\end{Def}

\begin{prop}
    If $\mathscr{A}$ is a dense $*$-subalgebra of $A$ on which $K^*(A)$ is $p$-summable, then there exists $\mathscr{A} \sub \mathcal{A} \sub A$ where, for each $n,\: M_n(\mathscr{A})$ is closed under the holomorphic functional calculus in $M_n(A)$ so that $\iota_*:K_*(\mathscr{A}) \to K_*(A)$ is an isomorphism, and on which $K^*(A)$ is $p$-summable \cite{C}. 
\end{prop}
Thus, when a $C^*$-algebra $A$ has uniformly $p$-summable $K$-homology, any index pairing between $K_*(A)$ and $K^*(A)$ can be computed via Connes' trace formulas for $p$-summable Fredholm modules because $\bigcup_{k \in \N} M_k(\mathscr{A})$ contains a representative of each $K_*(A)$ class.

\subsection{$KK$-Theory and the UCT}
The tools of Kasparov's $KK$-theory will be useful to translate results between $K$-theory and $K$-homology. For a full exposition on $KK$-theory, see \cite{K}.  Note that, for a separable $C^*$-algebra $A,\: KK_*(\C,A) \cong K_*(A)$ and $KK_*(A,\C) \cong K^*(A).$ $KK$-theory has an associative product. We utilize a case of this product called the \textit{cap product.}

\begin{Def}
    Let $A,B,D$ be separable $C^*$-algebras. There is an associative product $- \otimes_B - : KK(A,B) \times KK(B,D) \to KK(A,D)$ called the \textbf{cap product}.
\end{Def}

Additionally, the index paring is an example of the Kasparov product:
\begin{Ex}\label{IndexKKprod}
    Let $x \in K_*(A), y\in K^*(A).$ Then $\langle x,y \rangle_A = x \otimes_A y \in KK(\C,\C).$
\end{Ex}
\begin{Ex}\label{KKhom}
    Classes in $KK(A,B)$ represent `generalized morphisms' from $A$ to $B.$ When $\vp:A \to B$ is a $*$-homomorphism, $[\vp] = [B, \vp, 0]$ is a $KK_0(A,B)$ class. Further, the cap product satisfies: $[\vp] \otimes_B - = \vp^*(-) : KK(B, D) \to KK(A,D),$ while $- \otimes_A [\vp] = \vp_*(-):KK(D,A) \to KK(D,B),$ which we utilize in Section 6.
\end{Ex} 
\subsection{The Universal Coefficient Theorem} The tool we use to obtain results about $K$-homology from results on $K$-theory is the Universal Coefficient Theorem (UCT): \begin{Def}\label{UCT} We say that a $C^*$-algebra $A$ \textbf{satisfies the UCT} if, for all separable $C^*$-algebras $B,$ the following sequence is exact: \[0 \xrightarrow{} \text{Ext}_\Z^1(K_*(A), K_*(B)) \xrightarrow{\de} KK_*(A,B) \xrightarrow{\gm} \text{Hom}(K_*(A), K_*(B)) \xrightarrow{} 0.\] Recall that $\de$ has degree 1 and $\gm$ has degree 0.
\end{Def}
All $C^*$-algebras in this paper will satisfy the UCT, in particular, AF-algebras do, commutative algebras do, and the property of satisfying the UCT is preserved when taking crossed products by $\Z$ \cite{RS}.
\begin{Def}
    The map $\gm$ is given by $\gm(x):= - \otimes_A [x]:KK(\C,A) \to KK(\C, B)$ where $x \in KK(A,B).$
\end{Def} 
Thus, in the case of Example \ref{IndexKKprod}, we obtain: 
\begin{prop}
    If $A$ satisfies the UCT, $x \in K_*(A),$ and $y \in K^*(A),$ then $\langle x, y \rangle = \gm(y)(x).$
\end{prop}
\section{Cantor Minimal Systems}
The $C^*$-algebras we study in this paper arise from Cantor minimal systems. These are dynamical systems obtained from \textit{minimal homeomorphisms} on the \textit{Cantor set}.
    \begin{Def}
    A topological space $X$ is called a \textbf{Cantor set} if it is a totally disconnected, compact, metrizable space with no isolated points. (There is a unique space with these properties, so we refer to \textbf{the} Cantor Set.)
\end{Def}
We often use a symbolic representation of the Cantor set. \begin{Ex}\label{Cantoralph}
    Common models of the Cantor set include $\Omega^\Z, \Omega^\N,$ or a closed subset of either for a finite set $\Omega.$ 
    \begin{itemize}
        \item The metric on $\Omega^{\Z}$ (or a closed subset of $\Omega^\Z)$ is given by $d((x_i),(y_i)) = 2^{- \min\{|j|||x_j \neq y_j\}}.$
        \item The metric on $\Omega^\N$ (or a closed subset of $\Omega^\N$) is given by $d((x_i),(y_i)) = 2^{- \min\{j||x_j \neq y_j\} + 1}.$
    \end{itemize}  
\end{Ex}
This topology can also be given by a basis of clopen sets. 
\begin{Def}
    Suppose $X \sub \Omega^\Z$ is a Cantor set. Let $(x_n)_{n\in \Z} \in X$ and $a < b \in \Z.$ We denote by $x_{[a,b]} = x_a x_{a+1} \cdots x_{b}.$ Further $x_{(-\infty, b]} = \ldots x_{b-2} x_{b-1} x_b$ and $x_{[b,\infty)} = x_b x_{b+1} \ldots$ will denote the left and right tails of $x$ respectively.
\end{Def}
\begin{Def}
    For any $x \in X,\: a < b \in \Z, \: x_{[a,b]}$ will denote a \textbf{word} in X. For such a word $ \mu = x_{[a,b]}, \: |\mu| = b-a+1,$ is called the \textbf{length} of $\mu.$ The \textit{empty word}, denoted $\ee,$ is the unique word of length 0.
\end{Def}
\begin{Def}
    We denote by $Y$ the set of all words in $X: Y = \{x_{[a,b]} |x \in X, a< b \in \Z\} \cup \{\ee\}.$
\end{Def}
\begin{Def}
    If $\mu, \nu \in Y$ and $\nu = \mu a_1\cdots a_k$ for $a_1, \ldots, a_k \in \Omega,$ we say that $\mu$ is a \textbf{subword} of $\nu.$ Notice that the empty word is a subword of every word in $Y.$
\end{Def}
\begin{Def}\label{cylinder}
    For $\mu \in Y,$ denote by $C_\mu$ the \textbf{cylinder set defined by }$\mu$: \begin{itemize} \item If $X \sub \Omega^\N,$ $C_\mu = \{x \in X| x_{[0,|\mu| -1]} = \mu \}.$
    \item If $X \sub \Omega^\Z, |\mu|$ odd: $C_\mu = \{x \in X| x_{[-\frac{|\mu| -1}{2},\frac{|\mu| -1}{2}]} = \mu\}.$
    \item If $X \sub \Omega^\Z, |\mu|$ even: $C_\mu = \{x \in X| x_{[-\frac{|\mu|-1}{2},\frac{|\mu|}{2}]} = \mu\}.$
    \end{itemize}
\end{Def}
We collect a few facts below:
\begin{prop}
Let $\mu \in Y.$
    \begin{itemize}
        \item In the metric topology, each $C_\mu$ is clopen.
        \item  For any closed $X \sub \Omega^\Z$ (or $\Omega^\N$), the collection $\{ C_\mu | \mu \in Y\}$ forms a basis for the topology on $X.$
        \item Further, for each $n,$ $\mathscr{C}_n =\{C_\mu| \:|\mu| =n\}$ is a partition of $X.$ 
        \item These partitions are \textit{increasing}, in the sense that, for each $m > n,$ and each $Z \in \mathscr{C}_m,$ there is a $Z' \in \mathscr{C}_n$ such that $Z \sub Z'.$ 
    \end{itemize}
\end{prop}
\begin{prop}
    If $X$ is the Cantor set, \[K_0(C(X)) \cong C(X,\Z) \cong \text{span}\{\chi_E | E \text{ clopen } \sub X\}, \: \: K_1(C(X)) \cong \{0\}.\] Thus, by the UCT \[K^0 (C(X)) \cong \text{Hom}(C(X, \Z), \Z), \: \: K_1(C(X)) \cong 0,\] for example, see \cite{P} or \cite{B}.
\end{prop} 
\subsection{Odometers}
    Odometers are a particularly nice class of Cantor minimal systems. We can associate an odometer to each sequence of integers, $\{d_i\}_{i=1}^\infty, d_i \geq 2$ \cite{P}. 
\begin{Def}\label{Odometer}
    We let $X= \prod_{i=1}^\infty X_i$ where $X_i = \{0,1,\ldots, d_i -1\}$ and define $\vp: X \to X $ as follows. First, $\vp(d_1-1,d_2-1, \ldots, d_n -1, \ldots) = (0,0,0, \ldots).$ Now, suppose $(x_1,x_2, \ldots) \in X$ where for some $j, x_j \neq d_j-1.$ Then define $k = \min\{j| x_j \neq d_j -1\}$ and $\vp(x_1,x_2,\ldots) = (0,\ldots, 0, x_k +1, x_{k+1},x_{k+2}, \ldots),$ i.e. $\vp$ acts on $X$ via ``add $(1,0,0,\ldots)$ with carry.''
\end{Def}
\begin{Def}
    The topology on $X$ is induced by a metric $d:X\times X \to \R_{\geq 0}.$ The metric is defined by $d((x_i), (y_i)) = 2^{-\min\{j| x_j \neq y_j\} + 1}.$
\end{Def}
This topology makes $X$ into a Cantor set and makes $\vp$ a minimal homeomorphism. 
\subsection{Ordered Bratteli diagrams and Bratteli-Vershik Systems} 
Bratteli diagrams (and their associated Bratteli-Vershik systems) are a valuable tool for studying Cantor minimal systems. For further details on such systems see \cite{HPS}, \cite{PB}. In particular, every Cantor minimal system is conjugate to a Bratteli-Vershik system associated to a \textit{properly ordered Bratteli diagram} \cite[Theorem 4.6]{HPS}. \begin{Def}
    A \textbf{Bratteli diagram} $(V,E)$ consists of a set of vertices $V = \bigsqcup_{n=0}^\infty V_n,$ a set of edges $E = \bigsqcup_{n=1}^\infty E_n,$ and maps $s:E_n \to V_{n-1}, r: E_n \to V_n,$ called the \textbf{range} and \textbf{source} maps. Moreover, $V_n$ and $E_n$ are finite non-empty disjoint sets and $V_0 = \{v_0\}$ is a one-point set.
\end{Def}
Each edge $e_n \in E_n$ connects the vertex $s(e_n) \in V_{n-1}$ to the vertex $r(e_n) \in V_n.$ We assume that $s^{-1}(v) \neq \es$ for all $v \in V$ and $r^{-1}(v) \neq \es$ for all $v \not \in V_0.$ 
\begin{Def}
    If $|V_{n-1}| = k_{n-1}, |V_n| = k_n,$ then $E_n$ can be described by a $k_n \times k_{n-1}$ \textbf{transition matrix}, $S_n = [S_{ij}^n],$ where $S_{ij}^n$ is the number of edges connecting $v_i^n \in V_n$ with $v_j^{n-1} \in V_{n-1}.$
\end{Def} 

\begin{Def}\label{Bratpath}
    The \textbf{infinite path space} associated to the Bratteli diagram $(V,E)$ is \[X_{(V,E)} = \{(e_1,e_2, \ldots)|\: e_i \in E_i, r(e_i) = s(e_{i+1}) \text{ for all } i \geq 1\}.\]  We topologize $X_{(V,E)}$ using the subspace topology inherited from $\prod_{n=1}^\infty E_n,$ which is endowed with the product topology.
\end{Def}
\begin{Def}\label{Bratorder}
      An \textbf{ordered Bratteli diagram} $(V,E,\geq)$ is a Bratteli diagram $(V,E)$ together with a partial order $\geq$ in $E$ so that edges $e,e' \in E$ are comparable if and only if $r(e) = r(e').$
  \end{Def}  
  We denote by $E_{min},E_{max}$ the minimal and maximal edges of the poset $E.$

\begin{Def}\label{BV}
The \textbf{Vershik map} associated to a \textit{properly ordered Bratteli diagram} (see \cite{HPS},\cite{PB}) $(V,E,\geq)$ is the map, $T,$ given by $T(x_{max}) = x_{min}$ and, if $x = (e_1,e_2,\ldots) \neq x_{max}, k$ is the minimal number so that $e_k \not \in E_{max}, f_k$ is the successor of $e_k,$ and $(f_1,f_2, \ldots, f_{k-1})$ be the unique minimum path in $E_1 \circ E_2 \circ \cdots \circ E_{k-1}$ from $s(F_k)$ to $V_0.$ Then $T(e_1,e_2, \ldots) = (f_1,f_2, \ldots, f_k, e_{k+1}, e_{k+2}, \ldots).$ Then $(X_{(V,E)},T)$ is called the \textit{Bratteli-Vershik System} associated to $(V,E, \geq).$ \end{Def} 

Given a Cantor minimal system, $(X,\vp),$ \cite[Section 4]{HPS} shows how to construct a Bratteli-Vershik system conjugate to it. This is the same construction as used to produce the orbit-breaking subalgebras as in Section \ref{subalg}. We will use the following example of a Bratteli-Vershik system as our running example of a non-odometer Cantor minimal system:
\begin{Ex}\label{Bratexamp}
    Let $(X,\vp)$ be the Bratteli-Vershik system associated to the stationary Bratteli diagram $(V,E)$ where $|V_n| = 2$ for each $n \geq 1$ with transition matrix $S = \begin{bmatrix} 1 & 1 \\ 1 & 0 \end{bmatrix}.$ 
Following \cite{HPS} and \cite{B}, we label edges as follows:
 \begin{figure}[H] 
 \setlength{\belowcaptionskip}{-10pt}
    \[
\begin{tikzpicture}[node distance=3cm, 
  every node/.style={circle, draw, fill=black, inner sep=1.5pt}]

  \node (v6) at (-3,1) {};
  \node (v0) at (0,0) {};
  \node (v1) at (0,2) {};
  \node (v2) at (3,0) {};
  \node (v3) at (3,2) {};
  \node (v4) at (6,0) {};
  \node (v5) at (6,2) {};
  \node[draw=none,fill=none] (ellipsis) at (7,1) {\dots};
  \draw[-] (v6) to (v1);
  \draw[-] (v6) to (v0);
  \draw[-] (v1) to node[draw=none, fill=none,pos=0.4, above] {$1$} (v2);
  \draw[-] (v0) to node[draw=none, fill=none, midway, below] {$0$} (v2);
  \draw[-] (v0) to node[draw=none, fill=none, pos=0.4,below]   {$2$} (v3);
   \draw[-] (v3) to node[draw=none, fill=none, pos=0.4, above] {$1$} (v4);
  \draw[-] (v2) to node[draw=none, fill=none, midway, below] {$0$} (v4);
  \draw[-] (v2) to node[draw=none, fill=none, pos=0.4,below]   {$2$} (v5);
\end{tikzpicture}
\] 
\caption{Stationary Bratteli diagram with transition matrix $S = \begin{bmatrix} 1 & 1\\ 1 & 0 \end{bmatrix}$}
\end{figure}
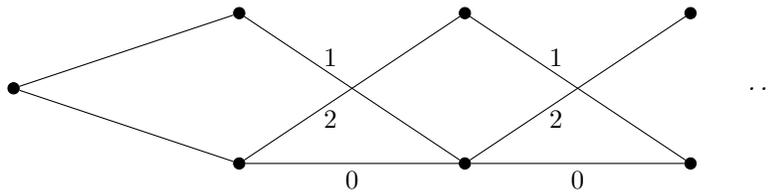
\end{Ex}
\subsection{Orbit-breaking Subalgebras and Embedding}
Our results on index pairings rely heavily on the orbit-breaking AF-algebras constructed by Putnam and the sequences of $C^*$-algebras and $K$-theory relating them to $C(X) \rtimes_\vp \Z.$ In this section, we describe these algebras and sequences. 
\subsection{``Orbit-breaking'' AF-Subalgebras}\label{subalg}
In \cite{P}, Putnam constructs an AF-algebra $A_{\{y\}} \sub C(X) \rtimes_\vp \Z$ obtained by ``breaking orbits at $\{y\}$.'' The algebra $A_{\{y\}}$ is the $C^*$-subalgebra of $C(X) \rtimes_\vp \Z$ generated by $C(X)$ and $uC_0(X-\{y\}).$ Specifically, $A_{\{y\}} \cong \br{\bigcup_n A_n}$ where $A_n \cong C^*(u\chi_{X-Z_n}, \{\chi_E | E \in P_n\})$ is finite dimensional, $Z_n$ is a decreasing sequence of clopen sets whose intersection is $\{y\},$ and $\{P_n\}$ is an increasing sequence of partitions of $X.$ The algebra satisfies: 
\begin{prop}{\cite[Theorem 4.1]{P}}\label{inclusK_0}
    Suppose $y \in X.$ Then the inclusion $j:A_{\{y\}} \to C(X) \rtimes_\vp \Z$ induces an isomorphism $j_*:K_0(A_{\{y\}}) \to K_0(C(X) \rtimes_\vp \Z).$
\end{prop}

\subsection{Embedding Cantor Minimal Systems in AF-algebras}\label{Embed}
In addition, $C(X) \rtimes_\vp \Z$ embeds back into $A_{\{y\}}$ unitally \cite{P}. \begin{prop}{\cite[Chapter 6]{P}}\label{Putnamembed}
    There is an embedding $\iota:C(X) \rtimes_\vp \Z \to A_{\{y\}}.$
\end{prop}In particular, this construction ensures that:
\begin{prop}{\cite[Theorem 6.7]{P}}\label{Putnamiso}
The map $\iota: C(X) \rtimes_\vp \Z \to A_{\{y\}}$ induces an isomorphism $\iota_*:K_0(C(X) \rtimes_\vp \Z) \to K_0(A_{\{y\}}).$ 
\end{prop}
Further, in the proof of Proposition \ref{Putnamiso}, Putnam shows the following: 
\begin{prop}\label{compidentity}
    With $j:A_{\{y\}} \to C(X) \rtimes_\vp \Z$ the inclusion and $\iota:C(X) \rtimes_\vp \Z \to A_{\{y\}}$ as in Proposition \ref{Putnamembed}, then $\iota_* \circ j_*:K_0(A_{\{y\}}) \to K_0(A_{\{y\}})$ is the identity.
\end{prop}
These $K$-theory results are of crucial importance for us. We also utilize the following to show how projections generating $K_0$ behave under the embedding:
\begin{prop}
    In the AF-structure described above, if $f \in C(X) \cap A_n \sub A_{\{y\}} \sub C(X) \rtimes_\vp \Z$ then $\iota(f) \in A_{n+1}.$
\end{prop}

For specifics on these embeddings see \cite{P}. We include an example, using the notation of \cite[Chapter 6]{P}:
\begin{Ex}\label{Embedexamp}
    Let $(X,\vp)$ be the Bratteli-Vershik system associated to the stationary Bratteli diagram $(V,E)$ where $|V_n| = 2$ for each $n \geq 1$ with transition matrix $S = \begin{bmatrix} 1 & 1 \\ 1 & 0 \end{bmatrix}.$ 
Let $y = \br{0}.$ $ Y_n = C_{0^n}$ and $P_n = \{C_{\mu_1}, \ldots, C_{\mu_l}\}$ where $\mu_1, \ldots, \mu_l$ are the paths of length $n+1$ starting from the initial vertex. Then we let $Z_n = C_{0^{n+1}},$ so that $P'_n = P_n.$ Then $\lambda(Z_1) = \{3,5\}$ so that $A_1 \cong M_5(\C) \oplus M_3(\C).$ Next, $\lambda(Z_2) = \{8,5\}$ so that $A_2 \cong M_8(\C) \oplus M_5(\C).$ In general, we have that, if $A_n \cong M_{n_1}(\C) \oplus M_{n_2}(\C), \: A_{n+1} = M_{n_1 + n_2}(\C) \oplus M_{n_1}(\C).$ The embedding of $A_n$ into $A_{n+1}$ is given by \[(T_1,T_2) \mapsto \left (\begin{bmatrix} T_1 & 0 \\ 0 & T_2 \end{bmatrix}, T_1 \right ).\] Thus, we have that \[v_1 = \left(\begin{bmatrix} 0 & 0 & 0 &0 & 1 \\ 1 & 0 & 0 & 0 &0 \\ 0 & 1 & 0 & 0 &0 \\ 0 & 0 & 1 & 0 & 0 \\ 0 & 0 & 0 & 1 & 0 \end{bmatrix}, \begin{bmatrix} 0 & 0 & 1 \\ 1 & 0 & 0 \\ 0 & 1 & 0 \end{bmatrix}\right ) \text{ while } v_2 = \left(\begin{bmatrix}  0 &0 & 0 & 0 & 0 & 0 & 0 & 1 \\ 1 & 0 & 0 & 0 & 0 &0 & 0 &0 \\ 0 & 1 & 0 & 0 & 0 & 0 & 0 & 0\\ 0 & 0 & 1 & 0 & 0 & 0 & 0 & 0\\ 0 & 0 & 0 &1 & 0 & 0 & 0 & 0 \\ 0 & 0 & 0 & 0 & 1 & 0 & 0 & 0\\ 0 & 0 & 0 & 0 & 0 & 1 & 0 & 0 \\ 0 & 0 & 0 & 0 & 0 & 0 & 1 & 0  \end{bmatrix}, \begin{bmatrix} 0 & 0 & 0 &0 & 1 \\ 1 & 0 & 0 & 0 &0 \\ 0 & 1 & 0 & 0 &0 \\ 0 & 0 & 1 & 0 & 0 \\ 0 & 0 & 0 & 1 & 0 \end{bmatrix}\right ).\] Then $v_2v_1^* \in A_2$ is given by \[v_2v_1^* = \left(e_{16} + e_{61} + \sum_{i= 2, i \neq 6}^8 e_{ii},\text{Id}_5\right).\] Thus, we can take \[z = \left (\frac{e^{i\frac{\pi}{4}}}{\sqrt{2}}(e_{11} + e_{66}) + \frac{e^{-i\frac{\pi}{4}}}{\sqrt{2}}(e_{16} + e_{61}) + \sum_{i= 2, i \neq 6}^8 e_{ii},\text{Id}_5 \right).\] Hence, $w_1 =z^2uzu^{-1}$ is given by \[w_1 = z^2 u z u^{-1} =  v_2v_1^* u z u^{-1} = \left(\begin{bmatrix} 0 & 0 & 0 & 0 & 0 & 1 & 0 & 0 \\ 0 & \frac{e^{i \frac{\pi}{4}}}{\sqrt{2}} & 0 & 0 & 0 & 0 & \frac{e^{ - i \frac{\pi}{4}}}{\sqrt{2}} & 0 & \\ 0 & 0 & 1 & 0 & 0 & 0 & 0 & 0 \\ 0 & 0 & 0 & 1 & 0 & 0 & 0 & 0 \\ 0 & 0 & 0 & 0 & 1 & 0 & 0 & 0 \\ 1 & 0 & 0 & 0 & 0 & 0 &0 & 0 \\ 0 & \frac{e^{-i \frac{\pi}{4}}}{\sqrt{2}} & 0 & 0 & 0 & 0 & \frac{e^{i \frac{\pi}{4}}}{\sqrt{2}} & 0  \\ 0 & 0 & 0 & 0 & 0 & 0 & 0 & 1\end{bmatrix}, \text{Id}_5 \right).\] Proceeding in this way, we obtain that: if $A_n = M_{n_1}(\C) \oplus M_{n_2}(\C)$ then \[z = \left (a e_{11} + b e_{1 (n_2 + 1)} + c e_{(n_2+1) 1} +d e_{(n_2+1)(n_2 +1)} + \sum_{i =2, \neq n_2 +1}^{n_1} e_{ii}, \text{Id}_{n_2}\right)\] where $a,b,c,d$ are chosen such that $z^{2^n}_{11} = z^{2^n}_{(n_2 + 1)(n_2 + 1)} = 0$ and $z^{2^n}_{1 (n_2 + 1)} = z^{2^n}_{(n_2 + 1) 1} = 1.$ Then we get that \[w_n = \left(v_{n+1}v_n^*\sigma(z^{2^n-1}) \sigma^2(z^{2^n-2}) \cdots, \text{Id}_{n_2}\right)\] where, for the $k \times k$ matrix $S:$ $\sigma(S)_{ij} = \begin{cases} S_{(i-1)(j-1)} & i,j >1 \\ S_{(i-1)k} & i >1, j =1 \\ S_{k(j-1)} & i=1, j > 1 \\ S_{kk} & i=j=1 \end{cases}.$ 
\end{Ex}

\section{Index Computations and Embeddings}\label{chapter:examples}
In this section, we show how, in the absence of a dense subalgebra on which we have finite summability, we can use embeddings into $C^*$-algebras with finitely summable $K$-homology to compute index pairings. We apply this to the case of Cantor minimal systems and their embeddings into AF orbit-breaking algebras in the next section. We also show we can exhaust any index pairing so long as we have enough cycles to surject onto the Hom-term in the UCT. We apply this result to Cantor minimal systems in the odd case in the next section.
\subsection{Embedding Results on $K$-homology}
In this section, we prove the following result.
 \begin{theorem}\label{maineventhm}
     Let $A$ and $B$ be $C^*$-algebras in the UCT class (Definition \ref{UCT}) with $K^*(A)$ uniformly $p$-summable on $\mathscr{A}.$ If there are $*$-homomorphims $\vp: A \to B, \psi: B \to A$ that satisfy: 
\begin{enumerate}
    \item $\vp$ and $\psi$ are injective,
    \item $\vp_*:K_0(A) \to K_0(B), \psi_*:K_0(B)\to K_0(A)$ are isomorphisms,
    \item $ \psi_* \circ \vp_* = \text{id}_{K_0(A)},$
    \item $K_{1}(A), K_{1}(B)$ are free abelian, and
    \item elements in $\psi^{-1}(\mathscr{A})$ generate $K_0(B),$
\end{enumerate}
then index pairings between $K_0(B)$ and $K^0(B)$ can be computed using Connes' trace formulas for $p$-summable cycles. 
 \end{theorem}  
The analogous result holds in the odd case (replacing all the $0$'s in the above theorem with $1$'s and vice-versa). For each of the lemmas of this section the analogous result holds in the odd case as well. We will apply the even case for Cantor minimal systems, so we focus on that. 

While partial conclusions can be made in the absence of one or more conditions, we rely on all of them to attain the main result. In the next section, we apply this to Cantor minimal systems. We begin as follows:
 
\begin{lemma}\label{Khomiso}
    Let $A$ and $B$ be $C^*$-algebras in the UCT class. Suppose $\vp:A \to B$ is a $*$-homomorphism such that $\vp_*:K_0(A) \to K_0(B)$ is an isomorphism. If $K_{1}(A), \: K_{1}(B)$ are free abelian, then $\vp^*:K^0(B) \to K^0(A)$ is an isomorphism. 
\end{lemma}
\begin{proof}
  Suppose $K_{1}(A), K_{1}(B)$ are free abelian. Then $\text{Ext}_\Z^1(K_{1}(A), K_0(\C)) \cong \text{Ext}_\Z^1(K_0(A), K_1(\C))  \cong \{0\}$ and the same for $B.$ The UCT (Definition \ref{UCT}) gives that \[0 \xrightarrow{} \text{Ext}_\Z^1(K_{1}(A), \Z) \xrightarrow{\de_A} KK_0(A,\C) \xrightarrow{\gm_A} \text{Hom}(K_0(A),\Z) \xrightarrow{} 0\] is exact and similarly for $B.$ Thus, $\gm_A:KK_0(A,\C) \to \text{Hom}(K_0(A),\Z)$ is an isomorphism, as is $\gm_B:KK_0(B,\C) \to \text{Hom}(K_0(B),\Z)$. Since $\vp_*:K_0(A) \to K_0(B)$ is an isomorphism, we have that $f \mapsto f \circ \vp_*: \text{Hom}(K_0(B),\Z) \to \text{Hom}(K_0(A),\Z)$ is an isomorphism as well. Using the naturality of the UCT \cite{B}, we obtain that the following diagram commutes: 
  \[\begin{tikzcd}
0 \arrow[r] & 0  \arrow[r] \arrow[d] & KK_0(B, \C)  \arrow[r, "\gm_B"] \arrow[d, "\vp^*"] & \text{Hom }(K_0(B), \Z) \arrow[r] \arrow[d, "f \mapsto f \circ \vp_*"] & 0 \\
0 \arrow[r] & 0 \arrow[r]   & KK_0(A, \C) \arrow[r, "\gm_A"]   & \text{Hom}(K_0(A), \Z) \arrow[r]   & 0
\end{tikzcd}.\]
Since $\gm_A, \gm_B, f \mapsto f \circ \vp_*$ are isomorphisms, $\vp^*:K^0(B) \to K^0(A)$ is as well.
\end{proof}
Note that we can remove the assumption of free abelian and still conclude: 
\begin{Cor}
    Let $A, B$ be $C^*$-algebras in the UCT class. Suppose $\vp:A \to B$ is a $*$-homomorphism such that $\vp_*:K_0(A) \to K_0(B)$ is an isomorphism. Then $\vp^*|_{\gm(K^0(B))}:\gm_B(K^0(B)) \to \gm_A(K^0(A)) \sub \text{Hom}(K_0(A),\Z)$ is an isomorphism. 
\end{Cor}
We return to the case where $K_{1}(A), K_{1}(B)$ are free abelian. Suppose, additionally, we have a map $\psi: B \to A$ such that $\psi_* \circ \vp_*:K_0(A) \to K_0(A) = \text{id}_{K_0(A)}$ (which implies $\vp_*,\psi_*$ are isomorphisms). Then we attain:
\begin{lemma}\label{Khomidentity}
    Let $A,B$ by $C^*$-algebras in the UCT class with $A \xrightarrow{\vp} B \xrightarrow{\psi} A$ such that $\vp^*:K^0(B) \to K^0(A)$ is an isomorphism, $K_{1}(A), K_{1}(B)$ are free abelian, and $\psi_* \circ \vp_* = \text{id}_{K_0(A)}:K_0(A) \to K_0(A).$ Then: 
    \begin{enumerate}
        \item $\vp_* \circ \psi_* = \text{id}_{K_0(B)},$
        \item $\psi^* \circ \vp^*= \text{id}_{K^0(B)},$ and
        \item $\vp^* \circ \psi^* = \text{id}_{K^0(A)}.$
    \end{enumerate} 
\end{lemma}
\begin{proof}
    For (1), observe that, since $K_0(A) \cong K_0(B),$ and $\vp_*$ is the right inverse to $\psi_*,$ it is also the left inverse. Specifically, for $x \in K_0(A):$ \[x = (\psi_* \circ \vp_*) \circ (\psi_* \circ \vp_*)(x)  = (\psi_* \circ (\vp_* \circ \psi_*) \circ \vp_*(x)).\]  Since $\psi_* \circ \vp_*$ is the identity on $K_0(A)$ and $\vp_*, \psi_*$ are isomorphisms, it must be that $\vp_* \circ \psi_*(\vp_*(x)) = \vp_*(x)$ so that $\vp_* \circ \psi_* = \text{id}_{K_0(B)}.$ Thus $\vp_* \circ \psi_* = \text{id}_{K_0(B)}.$ For (2), let $x \in K_0(B), \: y \in K^0(B)$ and $[\vp], [\psi]$ the classes $[B, \vp, 0] \in KK_0(A,B), \: [A,\psi, 0] \in KK_0(B,A),$  respectively. Using that $\vp_* \circ \psi_*$ is the identity on $K_0(B),$ Example \ref{KKhom}, and that the Kasparov product is associative, we have that: 
    \begin{align*}
        \gm_B(\psi^*\circ \vp^*(y))(x) &= x \otimes_B ([\psi] \otimes_A [\vp] \otimes_B y) \\
        &= (x \otimes_B [\psi] \otimes_A [\vp]) \otimes_B y \\
        &= \psi_* \circ \vp_* (x) \otimes_B y\\
        &= x \otimes_B y \\
        &= \gm_B(y)(x).
    \end{align*}
Since $\gm_B:K^0(B) \to \text{Hom}(K_0(B),\Z)$ is an isomorphism and $\gm_B(\psi^* \circ \vp^*(y)) = \gm_B(y), \: \psi^* \circ \vp^* = \text{id}_{K^0(B)}.$ Since $\pi_* \circ \vp_* = \text{id}_{K_0(A)},$ this computation proves (3) as well.
\end{proof}
Given this result, we are able to use the embeddings $\vp, \psi$ to compute index pairings.
\begin{lemma}\label{indexlemma}
    Let $A$ and $B$ be $C^*$-algebras with $*$-homomorphisms $\vp: A \to B, \psi: B \to A.$  Then \[\langle \psi_*(-), \vp^*(-)\rangle_{A} = \langle -, \psi^* \circ \vp^*(-) \rangle_B: K_*(B) \times K^*(B) \to \Z.\]
\end{lemma}
\begin{proof}
    Let $x \in KK_*(\C,B), y \in KK^*(B,\C)$ and denote by $[\vp], [\psi]$ the classes $[B, \vp, 0] \in KK_0(A,B), \: [A,\psi, 0] \in KK_0(B,A)$ respectively. Observe that, using the fact that the index pairing coincides with the Kasparov product (Example \ref{IndexKKprod}), Example \ref{KKhom}, and the associativity of the Kasparov product: 
    \begin{align*} 
   \langle \psi_*(x), \vp^*(y) \rangle &= (x \otimes_B [\psi]) \otimes_A ([\vp] \otimes_B y) \\ 
   &= x \otimes_B ([\psi] \otimes_A [\vp] \otimes_B y) \\
   &= \langle x, \psi^* \circ \vp^*(y) \rangle \in KK_0(\C,\C)
    \end{align*}
\end{proof}
\begin{Cor}\label{indexresult}
    Let $A,B$ be $C^*$-algebras with $*$-homomorphisms, $\vp:A \to B, \psi: B \to A$ such that $\psi^* \circ \vp^* = \text{id}_{K^0(B)}.$ Then $\langle \vp_*(-), \psi^*(-)\rangle = \langle -, - \rangle: K_0(A) \times K^0(A) \to \Z.$
\end{Cor}
 \begin{proof}
     Let $x \in KK_0(\C,B), y \in KK^0(B,\C),$ and denote by $[\vp]$ and $ [\psi]$ the classes $[B, \vp, 0] \in KK_0(A,B)$ and $ [A,\psi, 0] \in KK_0(B,A),$  respectively. Observe that, using Lemma \ref{indexlemma}: \begin{align*}
        \langle \psi_*(x), \vp^*(y)\rangle &= \langle x, \psi^* \circ \vp^*(y) \rangle\\
        &= \langle x,y \rangle \in KK_0(\C,\C).
    \end{align*}
    
\end{proof}
\begin{Cor}\label{projectioncor}
Suppose that $A,B,\vp, \psi$ satisfy the assumptions of Lemma \ref{Khomiso}. Suppose further that $K^0(A)$ is uniformly $p$-summable on $\mathscr{A}$ and that $p \in B$ is a projection such that $\psi(p) \in \mathscr{A}.$ Then, for each $z \in K^0(A),$ there is an even Fredholm module $\left(H_+ \oplus H_-, \rho_+ \oplus \rho_-, \begin{bmatrix} 0 & U \\ V & 0 \end{bmatrix}\right)$ representing $z$ such that $\left \langle [p], \left[ H_+ \oplus H_-, \rho_+ \oplus \rho_-, \begin{bmatrix} 0 & U \\ V & 0 \end{bmatrix}\right ] \right \rangle = \left \langle \vp_*[p], \psi^*\left[ H_+ \oplus H_-, \rho_+ \oplus \rho_-, \begin{bmatrix} 0 & U \\ V & 0 \end{bmatrix}\right] \right\rangle$ can be computed using Connes' trace formula for even $p$-summable cycles. The analogous result holds in the odd case replacing projections with unitaries.
\end{Cor}  
\begin{proof}
    Since $K^0(A)$ is uniformly $p$-summable on $\mathscr{A},$ for each $z \in K^0(B), \: \vp^*(z)$ has a representative that is $p$-summable on $\mathscr{A}.$ While $\psi^{-1}(\mathscr{A})$ may or may not be dense in $B$, any projection $p$  in $\vp^{-1}(\mathscr{A})$ has a representative that has $L^p$ commutator with $\vp(p).$ Thus, we can compute the pairing of $p$ and $z$ which, by Lemma \ref{Khomiso} is equal to the pairing between $\vp_*(p)$ and $\psi^*(z),$ using Connes' trace formula for $p$-summable cycles \cite{C}. 
\end{proof}
The main theorem follows:
\begin{Cor}\label{indexsumcomp}
    If $A,B, \vp, \psi$ satisfy the assumptions of Theorem \ref{maineventhm}, then any index pairing between a class in $K^0(A)$ and one in $K_0(A)$ can be computed using Connes' trace formulas for $p$-summable cycles.
\end{Cor}
\begin{proof}
    Since elements of $\psi^{-1}(\mathscr{A})$ generate $K_0(B),$ for each class in $K_0(B)$ we can find an element that gets mapped into $\mathscr{A}$ by $\psi$ and thus apply the previous corollary.
\end{proof}
\subsection{Index Pairings and the UCT}
The next result shows that, for any $K$-homology class over a $C^*$-algebra that satisfies the UCT, one can compute index pairings with this class by computing with any other class that has the same image under $\gm:KK(A,\C) \to \text{Hom}(K_*(A), K_*(\C)).$
\begin{lemma}\label{IndexUCT}
    Let $A$ be in the UCT class. Let $x \in K^*(A).$ Then, for any $z \in K^*(A)$ such that $\gm(z) = \gm(x)$ and any $y \in K_*(A),$ \[\langle z, y \rangle = \langle x,y \rangle. \] 
\end{lemma}
\begin{proof}
    For $x \in K^*(A),y \in K_*(A)$ \[\gm(x)(y) = x \otimes_A y = \langle x, y \rangle.\] Thus, if $\gm(x) = \gm(z), \:\langle x,y \rangle = \gm(x)(y) = \gm(z)(y) = \langle z,y \rangle.$
\end{proof}
\begin{Cor}\label{IndexUCTCor}
    Let $x \in K^*(A).$ Suppose there is a $z \in K^*(A)$ and $\mathscr{A} \sub A$ such that $z$ is $p$-summable on $\A$ and $\gm(x) = \gm(z).$ Then pairings of $x$ with $K_*(A)$ classes represented by elements in $\A$ can be computed using one of Connes' trace formulas for $p$-summable cycles.
\end{Cor}
\begin{proof}
    Pairings with $z$ and $K_*(A)$ classes represented by elements in $\A$ can be computed using Connes' trace formula for $p$-summable cycles. By Lemma \ref{IndexUCT}, this pairing equals the pairing of this $K_*(A)$ class with $z.$
\end{proof}
\section{Index Pairing Results for Cantor Minimal Systems}
In this section, we apply the results of the previous section to Cantor minimal systems and their orbit-breaking subalgebras. We show:
\begin{theorem}\label{mainthm}
    Suppose $X$ is a Cantor set and $\vp: X \to X$ a minimal homeomorphism. Then any index pairing between an element in $K^*(C(X) \rtimes_\vp \Z)$ and one in $K_*(C(X) \rtimes_\vp \Z)$ can be computed using Connes' trace formulas for $p>0$ summable cycles.
\end{theorem}
The even and odd cases are handled separately, as in the previous section. Recall the set-ups of Sections \ref{subalg} and \ref{Embed}. We have that $A_{\{y\}},$ the subalgebra of $C(X) \rtimes_\vp \Z$ attained by breaking orbits at a point $y$, is AF, with structure described in Section \ref{subalg}. We also have that $C(X) \rtimes_\vp \Z$ can be embedded into $A_{\{y\}}$ via the map $\iota:C(X) \rtimes_\vp \Z \to A_{\{y\}}$ described in Section \ref{Embed}. Letting $j: A_{\{y\}} \to C(X) \rtimes_\vp \Z$ denote the inclusion, we have \[A_{\{y\}} \xrightarrow{j} C(X) \rtimes_\vp \Z \xrightarrow{\iota} A_{\{y\}}\] which induces \[K_0(A_{\{y\}}) \xrightarrow{j_*} K_0(C(X) \rtimes_\vp \Z) \xrightarrow{\iota_*} K_0(A_{\{y\}}).\] Then, Proposition \ref{inclusK_0} gives that $j_*$ is an isomorphism, 
 Proposition \ref{Putnamiso} gives that $\iota_*$ is an isomorphism, and Proposition \ref{compidentity} gives that $j_* \circ \iota_* = \text{id}_{K_0(A)}.$ Since $A_{\{y\}}$ is AF, we have that $K_1(A_{\{y\}}) \cong \{0\}.$ Further, \cite[Theorem 1.1]{P} gives that $K_1(C(X) \rtimes_\vp \Z) \cong \Z$ and is generated by $[u]$ where $u \in C(X) \rtimes_\vp \Z$ is the unitary implementing $\vp.$ Finally, since $A_{\{y\}}$ is AF, it has uniformly $p$-summable $K$-homology for all $p >0$ on the union of finite-dimensional subalgebras dense in it \cite[Section 4]{R}. The AF-filtration of $A_{\{y\}}$ is detailed in Section \ref{subalg}. These are the primary results we need to apply the results from the previous section.

In this section, we use the UCT to relate the $K$-homology of the orbit-breaking AF-algebras, $A_{\{y\}},$ to the $K$-homology of the crossed product algebras, $C(X) \rtimes_\vp \Z,$ for Cantor minimal systems. We use the results collected above to apply the results of the previous section to the Cantor minimal system setting. To begin:
 \begin{prop}
     Suppose $(X,\vp)$ is a Cantor minimal system and $A_{\{y\}}$ is the subalgebra of $C(X) \rtimes_\vp \Z$ obtained by breaking orbits at the point $y,$ see Section \ref{subalg}. Then $K^0(A_{\{y\}}) \cong K^0(C(X) \rtimes_\vp \Z)$  and $K^1(A_{\{y\}}) \oplus \Z \cong K^1(C(X) \rtimes_\vp \Z) .$
 \end{prop}
 \begin{proof}
We will use the UCT to prove this result. Observe that both $A_{\{y\}}$ and $C(X) \rtimes_\vp \Z$ satisfy the UCT. This is because $A_{\{y\}}$ is AF and $C(X) \rtimes_\vp \Z$ is a crossed product of a commutative algebra by $\Z$ \cite{B}. Per Proposition \ref{inclusK_0},
we have that $K_0(A_{\{y\}}) \cong K_0(C(X) \rtimes_\vp \Z).$ We also have that $K_1(A_{\{y\}}) \cong 0,$ since $A_{\{y\}}$ is AF. Additionally, $K_1(C(X) \rtimes_\vp \Z) \cong \Z$ \cite[Theorem 1.1(ii)]{P}. Then, the UCT gives that the following sequence is exact: \[ 
0 \xrightarrow{} \text{Ext}_\mathbb{Z}^1(K_*(C(X)\rtimes_\varphi \mathbb{Z}), K_{*+1}(\mathbb{C}))  \xrightarrow{\de} K^0(C(X) \rtimes_\varphi \mathbb{Z}) \xrightarrow{\gm} \text{Hom}(K_*(C(X) \rtimes_\varphi \mathbb{Z}), K_*(\mathbb{C}))  \xrightarrow{} 0.
\]
Since $\Z$ is free abelian, $K_0(\C) \cong \Z,$ and $K_1(\C) \cong 0,$ we have \[\text{Ext}_\mathbb{Z}^1(K_1(C(X)\rtimes_\varphi \mathbb{Z}), K_{0}(\mathbb\C)) \cong \text{Ext}_\mathbb{Z}^1(\Z,\Z) \cong 0.\] Further, \[K^0(C(X) \rtimes_\vp \Z) \cong \text{Hom}(K_0(C(X) \rtimes_\vp \Z), \Z) \cong \text{Hom}(K_0(A_{\{y\}}), \Z) \cong K^0(A_{\{y\}}).\] Plugging in for $K^1,$ we see that $K^1(A_{\{y\}}) \cong \text{Ext}_\Z^1(K_0(A_{\{y\}}), \Z),$ while \[0 \xrightarrow{} K^1(A_{\{y\}}) \xrightarrow{} K^1(C(X) \rtimes_\vp \Z) \xrightarrow{} \Z \xrightarrow{} 0\] is exact. Since this sequences splits (though not naturally), $K^1(C(X) \rtimes_\vp \Z) \cong K^1(A_{\{y\}}) \oplus \Z.$
 \end{proof}

\subsection{Even Case}
We seek to compute even index pairings for $C(X) \rtimes_\vp \Z$ using summable cycles. To do so, we use the orbit-breaking algebra $A_{\{y\}}.$ Because we are interested in the cycles representing classes in $K^0,$ we hope to identify a map inducing an isomorphism between $K^0(A_{\{y\}})$ and $K^0(C(X) \rtimes_\vp \Z).$ As we would hope, the inclusion has this property:
\begin{lemma}\label{incluiso}
    The map $\iota^*: K^0(A_{\{y\}}) \to K^0(C(X) \rtimes_\vp \Z)$ is an isomorphism.
\end{lemma}
\begin{proof}
       This is an application of Lemma \ref{Khomiso}, noting that $\iota_*:K_0(A_{\{y\}}) \to K_0(C(X) \rtimes_\vp \Z)$ is an isomorphism and $K_1(A_{\{y\}}) \cong \{0\}, K_1(C(X) \rtimes_\vp \Z) \cong \Z.$
\end{proof}
Thus, we can pull back $K$-homology classes along $\iota^*$ to exhaust $K^0(C(X) \rtimes_\vp \Z).$ We seek to compute index pairings using these cycles. Towards this end:
\begin{lemma}\label{khomid}
    Letting $j:A_{\{y\}} \to C(X) \rtimes_\vp \Z$ denote the inclusion and $\iota: C(X) \rtimes_\vp \Z \to A_{\{y\}}$ be given by Proposition \ref{Putnamembed}. Then $\iota^* \circ j^*: K^0(C(X) \rtimes_\vp \Z) \to K^0(C(X) \rtimes_\vp \Z)$ is the identity. 
\end{lemma}
\begin{proof}
   This is an application of Lemma \ref{Khomidentity}, noting that $K_1(A_{\{y\}}) \cong \{0\}, \: K_1(C(X) \rtimes_\vp \Z) \cong \Z,$ and $\iota_* \circ j_* = \text{id}_{K_0(A)}$ from \cite[proof of Theorem 6.7]{P}.
\end{proof}

Now, we can utilize this result to compute index pairings between the two algebras.
\begin{Cor}
    Let $x \in K_0(A_{\{y\}})$ and $z \in K^0(A_{\{y\}}).$ Then \[\langle \iota_*(x), j^*(y) \rangle_{C(X) \rtimes_\vp \Z} =  \langle x, y \rangle_{A_{\{y\}}}.\]
\end{Cor}
\begin{proof}
    This is an application of Lemma \ref{indexresult} with Lemma \ref{khomid}.
\end{proof}
\begin{Cor}\label{pairingcantor}
    Let $x \in K_0(C(X) \rtimes_\vp \Z)$ and $y \in K^0(C(X) \rtimes_\vp \Z).$ Then \[\langle \iota_*(x), j^*(y) \rangle_{A_{\{y\}}} = \left \langle x, y \right \rangle_{C(X) \rtimes_\vp \Z}.\]
\end{Cor}
\begin{proof}
    This is a direct application of Lemma \ref{indexresult} with Lemma \ref{khomid}.
\end{proof}
Recall that $K_0(C(X) \rtimes_\vp \Z)$ is generated by equivalence classes of projections in $C(X,\Z).$ Further, recall that $\iota_*$ and $j^*$ are isomorphisms in degree 0. Thus, we can compute index pairings of classes of projections with even Fredholm modules on $C(X) \rtimes_\vp \Z$ by pushing forward projections from $A_{\{y\}}$ and pulling back Fredholm modules to $A_{\{y\}}.$
Since $A_{\{y\}}$ is AF, its $K$-homology is uniformly finitely summable on the union of finite-dimensional algebras dense in it \cite[Section 4]{R}. We seek to use this to apply Corollary \ref{indexsumcomp} to compute index pairings over $C(X) \rtimes_\vp \Z.$ We use the AF-filtration, $A_{\{y\}} = \bigcup_n A_n$ where each $A_n = A(Z_n, P_n) \cong C^*(\{\chi_E| E \in P_n\}, u C_0(X-Z_n))$ for $Z_n$ clopen and $P_n$ a partition is finite-dimensional. Then we have:
\begin{lemma}\label{projections}
 Let $\iota: C(X) \rtimes_\vp \Z$ be as in Proposition \ref{Putnamembed}. Then $\iota^{-1}\left(\bigcup_n A_n\right) \cap C(X)$ is dense in $C(X).$ Further, projections in $\iota^{-1}\left(\bigcup_n A_n\right)$ generate $K_0(C(X) \rtimes_\vp \Z).$ 
\end{lemma} 
\begin{proof}
Observe that, for each $n, \: A_n$ contains $\chi_E$ for each $E \in P_n.$ By the construction of Proposition \ref{Putnamembed}, the sequence of partitions $\{P_n\}_n$ is increasing and its union generates the topology of $X$ \cite{P}. Then, by \cite[Lemma 6.4]{P}, if $f$ is a function in $\text{span}\{\chi_E| \: E \in P_n\},\: f$ commutes with $w_m$ for all $m > n,$ where the $w_m$ are unitaries in $A_{m+1}$ given in the construction of \cite[Chapter 6]{P}. Thus, \[\iota(f) = \lim_m (w_1^{-1} \cdots w_m^{-1} f w_m \cdots w_1) = w_1^{-1} \cdots w_n^{-1} f w_n\cdots w_1 \in A_{n+1},\] as $f \in A_{n+1}$ and $w_i \in A_{n+1},$ for each $1 \leq i \leq n.$ Thus $\iota\left(\bigcup_{n \in \N} \text{span}\{\chi_E |E \in P_n\}\right) \sub \bigcup_{n \in \N} A_n.$ Further, $\bigcup_{n \in \N} \text{span}\{\chi_E |\:E \in P_n\}$ is dense in $C(X).$ Now, recall that $K^0(C(X) \rtimes_\vp \Z) = \text{im} (j_*(C(X)) \cong \ker (1- \alp_*)(C(X))$ where $j:C(X) \to C(X) \rtimes_\vp \Z$ is the inclusion \cite{B}. Since $K_0(C(X)) \cong C(X, \Z),$ it is generated as a group by $\bigcup_{n \in \N} \text{span}\{\chi_E |\:E \in P_n\}.$ Thus, $K_0(C(X) \rtimes_\vp \Z)$ is generated by the classes of projections in $\iota^{-1}(\bigcup_{n \in \N} A_n).$ 
\end{proof}
Since $K^*(A_{\{y\}})$ is uniformly $p$-summable for all $p>0$ on $\iota^{-1}(\cup_n A_n)$, by Corollary \ref{indexsumcomp}, we obtain the following:
\begin{theorem}\label{eventhm}
    Suppose $(X,\vp)$ is a Cantor minimal system. Then any index pairing between a class in $K_0(C(X) \rtimes_\vp \Z)$ and one in $K^0(C(X) \rtimes_\vp \Z)$ can be computed using Connes' trace formula for even cycles with $p>0.$
\end{theorem}
Returning to the system of Example \ref{Bratexamp} and recalling the embedding of Example \ref{Embedexamp}:
\begin{Ex}
    By \cite{HPS}, the dimension group associated to $(V,E)$ (as in Example \ref{Bratexamp}) is \[K_0(V,E) \cong \varinjlim(\Z^2 \xrightarrow{S} \Z^2 \xrightarrow{S} \Z^2 \cdots) \cong \Z^2,\] since $S$ is invertible over $\Z.$ Then, by \cite[Theorem 5.4(2)]{HPS}, $K_0(C(X) \rtimes_\vp \Z) \cong \Z^2,$ and \cite[Theorem 1.1(i)]{P}, $K_1(C(X)\rtimes_\vp \Z) \cong \Z.$ Since $K_0,K_1$ are free abelian, the UCT gives that $K^0(C(X) \rtimes_\vp \Z) \cong K_0(C(X) \rtimes_\vp \Z) \cong \Z^2, \: K^1(C(X) \rtimes_\vp \Z) \cong K^1(C(X) \rtimes_\vp \Z) \cong \Z.$  Now, denote elements of $A_{\{y\}}$ (with AF-filtration as in Example \ref{Embedexamp}) by $[a,b,c]$ where $a \in M_{c_1}(\C), b \in M_{c_2}(\C), c \in \N.$ In this notation, $c$ indicates the level of the filtration. Then $K_0(C(X) \rtimes_\vp \Z)$ is generated by $\left [\begin{bmatrix} 1 & 0 & 0 & 0& 0\\ 0 & 0 & 0 & 0 & 0 \\ 0 & 0 & 0 & 0 & 0 \\ 0& 0 &0 &0 &0 \\ 0& 0 & 0 & 0 &0  \end{bmatrix}, 0, 1 \right ]$ and $\left [0, \begin{bmatrix} 1 & 0 & 0 \\ 0 & 0 & 0 \\ 0 & 0 & 0 \end{bmatrix}, 1 \right ].$ That is, $K_0$ is generated by the equivalence classes of the upper left projections in each coordinate in the 1st level of the filtration. Then $K^0(C(X) \rtimes_\vp \Z)$ is generated by the classes of the Fredholm modules which pair to $1$ with one generator and $0$ to the other. We can construct such generators using the technique of \cite[Theorem 4.2.1]{R}. For example, we can define Fredholm modules of the form $\left(H \oplus H, \vp_1 \oplus \vp_2, \begin{bmatrix} 0 & 1 \\ 1 & 0 \end{bmatrix}\right).$ Then $\iota^*\left(H \oplus H, \vp_1 \oplus \vp_2, \begin{bmatrix} 0 & 1 \\ 1 & 0 \end{bmatrix}\right)$ is $p$-summable for all $p>0$ on $\iota^{-1}\left(\bigcup_n A_n\right)$ and \[\left \langle j_*([p]) ,\iota^*\left(\left[\vp_1 \oplus \vp_2, \begin{bmatrix} 0 & 1 \\ 1 & 0 \end{bmatrix}\right]\right) \right \rangle = \left \langle [p], \left[\vp_1 \oplus \vp_2, \begin{bmatrix} 0 & 1 \\ 1 & 0 \end{bmatrix}\right] \right \rangle\] for each $[p] \in K_0(A_{\{y\}}).$

We can use these embeddings and Theorem \ref{maineventhm} to compute index pairings between elements of $K_0(C(X) \rtimes_\vp \Z)$ and $K^0(C(X) \rtimes_\vp \Z)$ using Connes' trace formulas.
\end{Ex}
\subsection{Odd Case}
We now seek to compute all index pairings with odd classes using summable cycles using Lemma \ref{IndexUCT}. 
 Our goal is to show that:  
\begin{lemma}\label{Cantorodd}
    Let $z \in \text{Hom}(K_1(C(X) \rtimes_\vp \Z), \Z).$ Let $\A = \text{span}\{\chi_E | E \text{ clopen } \sub X \}.$ Then there is a class $x \in K^1(C(X)\rtimes_\vp \Z)$ such that $\gm(x) = z$ and $x$ can be represented by a cycle that is $p$-summable for all $p >0$ on $C_c(\Z, \A).$ 
\end{lemma}
We have that:
\begin{Cor}\label{oddinclusion}
    The cokernel of $\iota^*:K^1(A_{\{y\}}) \to K^1(C(X) \rtimes_\vp \Z)$ is $\Z.$ 
\end{Cor}
\begin{proof}
    Recall that $\iota_*: K_0(C(X) \rtimes_\vp \Z) \to K_0(A_{\{y\}})$ is an isomorphism while $K_1(C(X) \rtimes_\vp \Z) \cong \Z$ and  $K_1(A_{\{y\}}) \cong 0.$ Since the UCT is natural in each variable \cite{B}, the following diagram commutes:  \[\begin{tikzcd}[sep=small]
0 \arrow[r] & \text{Ext}_\Z^1(K_0(C(X) \rtimes_\vp \Z), \Z)  \arrow[r, "\de"] & KK^1(C(X) \rtimes_\vp \Z, \C) \arrow[r] & \Hom(K_1(C(X) \rtimes_\vp \Z), \Z) \arrow[r]  & 0 \\
0 \arrow[r] & \text{Ext}_\Z^1(K_0(A_{\{y\}}), \Z) \arrow[r, "\cong"]  \arrow[u, "\cong"] & KK^1(A_{\{y\}}, \C) \arrow[r]  \arrow[u, "\iota^*"] & \Hom(K_1(A_{\{y\}}), \Z) \arrow[r] \arrow[u, "0"]  & 0
\end{tikzcd}.\] Thus $\iota^*$ is an isomorphism onto $\de(\text{Ext}_\Z^1(K_0(C(X) \rtimes_\vp \Z), \Z)) \sub K^1(C(X) \rtimes_\vp \Z),$ which has cokernel $\Z.$
\end{proof}
\begin{Cor}\label{oddquotient}
    Let $q: K^1(C(X) \rtimes_\vp \Z) \to K^1(C(X)\rtimes_\vp \Z)/ \iota^*(K^1(A_{\{y\}}))$ denote the quotient map. Let $x \in K^1(C(X) \rtimes_\vp \Z)$ and $u$ the unitary inducing $\vp.$ Then $q(x) \mapsto \langle [u], x \rangle$ under the identification $K^1(C(X)\rtimes_\vp \Z)/ \iota^*(K^1(A_{\{y\}})) \cong \Z.$
\end{Cor}
\begin{proof}
    Observe that  $K_1(C(X) \rtimes_\vp \Z) \cong \Z$ and is generated by $[u],$ \cite[Theorem 1.1(i)]{P}. Thus, if $x \in K^1(C(X)\rtimes_\vp \Z), \: \gm(x) \in \text{Hom}(K_1(C(X) \rtimes_\vp \Z),\Z)\cong \Z$ is determined by $\gm(x)([u]) := [u] \otimes_{C(X) \rtimes_\vp \Z} x =\langle [u], x \rangle,$ using the fact that the index pairing coincides with the Kasparov product. Since $\text{Hom}(K_1(C(X) \rtimes_\vp \Z),\Z) \cong K^1(C(X)\rtimes_\vp \Z)/ \iota^*(K^1(A_{\{y\}}))$ by the previous corollary, the result follows. 
\end{proof}
Then we have:
\begin{Cor}
   Let $x \in K^1(C(X) \rtimes_\vp \Z).$ Suppose $z \in K^1(C(X) \rtimes_\vp \Z)$ is such that $\gm(z) = \gm(x) \in \text{Hom}(K_1(C(X) \rtimes_\vp \Z), \Z).$ Then there is a $w \in \iota^*(K^1(A_{\{y\}}))$ such that $x = w \oplus z$
\end{Cor} 
\begin{proof}
    This follows from Corollary \ref{oddquotient} and the fact that the UCT exact sequence splits.
\end{proof} 
Motivated by these results, we construct odd finitely summable cycles whose images in the UCT exact sequence exhaust $\text{Hom}(K_1(C(X) \rtimes_\vp \Z), \Z)$ so that we can compute odd index pairings using Lemma \ref{IndexUCT} and Corollary \ref{IndexUCTCor}. We construct general nontrivial Fredholm modules on crossed products of commutative $C^*$-algebras. We start with a diagonal representation and use the covariant representation construction of Proposition \ref{hat}. We begin in the general case of a crossed product of a commutative $C^*$-algebra by $\Z.$ Then we apply this to the case of Cantor minimal systems. In the general setting, we have a commutative $C^*$-algebra, $A,$ and an automorphism $\vp: A \to A.$ We also have a representation $\pi:A \to B(H)$ on a separable Hilbert Space $H.$ Specifically, $\pi$ is such that there is a basis, $(\de_n)_{n \in J}$ for $J$ countable, with respect to which $\pi(a)$ is diagonal for each $a \in A.$ In this case, we are able to use Proposition \ref{hat} to produce odd Fredholm modules on $A \rtimes_\vp \Z$ where $A \rtimes_\vp \Z$ is represented on $\ell^2(\Z) \otimes H.$ These Fredholm modules will be finitely summable on the subalgebra \[C_c(\Z,A) := \left\{\sum_{k=L}^K a_ku^k \: \: L \leq K \in \Z, a_k \in A \right\} \text{ where } u \text{ implements } \vp.\] Denote by $e_m \in \ell^2(\Z)$ the vector $e_m(n) = \begin{cases} 1 & n=m \\ 0 & \text{ else}\end{cases}.$ Thus, $(e_m)_{m \in \Z}$ is an orthonormal basis for $\ell^2(\Z)$ so that $(e_m \otimes \de_n)_{(n,m) \in \Z^2}$ is an orthonormal basis for $\ell^2(\Z) \otimes H.$ Then we have:
\begin{prop}\label{CPodd}
Let $A$ be a commutative $C^*$-algebra and $\vp:A \to A$ an automorphism. Suppose that $\pi:A \to B(H)$ is a representation on the separable Hilbert space, $H,$ which satisfies that there is a countable basis $(\de_n)_{n \in J}$ for $H$ such that, for each $a \in A, \: n \in J, \: \pi(a)(\de_n) = \lambda_n \de_n$ for $\lambda_n \in \C.$ Next, let $N \sub J$ be a finite set and define $P_N \in B(\ell^2(\Z) \otimes H)$ to be projection onto the subspace spanned by those $e_m \otimes \de_n$ where $m > 0$ and $n \in N.$ Then $(\ell^2(\Z) \otimes H, \hat{\pi}, 2 P_N -1)$ is an odd Fredholm module over $A \rtimes_\vp \Z$ that is $p$-summable for all $p >0$ on $C_c(\Z, A).$ 
\end{prop}
\begin{proof}
 We show that, for $f =\ds \sum_{k=L}^K a_ku^k\in C_c(\Z,A), \: [\hat{\pi}(f),P_N]$ is finite rank. Observe that \[\hat{\pi}(a_ku^k) P_N(e_m \otimes \de_n) = \begin{cases} 0 & m\leq 0\text{ or } n \not \in N \\ e_{m+k} \otimes \pi(\vp^{-k}(a_k))\de_m & \text{ else } \end{cases}\] while \[P_N\hat{\pi}(a_ku^k)(e_m \otimes \de_n) = \begin{cases} 0 & m+k \leq 0 \text{ or } n \not \in N \\ e_{m+k} \otimes \pi(\vp^{-k}(a_k))\de_m &  \text{ else }\end{cases}.\] Thus, \[\left[\hat{\pi}\left(\sum_{k=L}^K a_ku^k\right),P_N\right](e_m \otimes \de_n)= \begin{cases} 0 & \text{ if }n \not \in N \text{ or } m \not \in [L,K] \\ \sum_{m >0,\: -m \geq k \geq L} e_{m+k} \otimes \pi(\vp^{-k}(a_k))\de_n & K \geq m >0, n \in \N \\ - \sum_{m \leq 0,\: K \geq k > |m|} e_{m+k} \otimes \pi(\vp^{-k}(a_k))\de_n & 0 \geq m \geq -K, n \in N \end{cases}.\] Because the representation $\pi$ is by diagonal operators with respect to $(e_n),$ for each $k$ and $n, \\ \pi(\vp^{-k}(a_k)) e_n = \lambda_{a_{k_n}} e_n$ for some $\lambda_{a_{k_n}} \in \C.$ Thus, for $f \in C_c(\Z,A),$ the operator $[\hat{\pi}(f), P_N]$ (and thus $[\hat{\pi}(f), 2P_N -1])$ has rank bounded by $(K-L+1)|N|.$ Now, following \cite[comment above Definition 2.3.11]{R}, we show that the commutator of the image of a general element of $A \rtimes_\vp \Z$ under $\hat{\pi}$ with $P_N$ is compact. Towards this, suppose $a \in A \rtimes \Z$ and $(a_n) \in C_c(\Z,A)^\N$ where $\lim_{n \to \infty} (a_n) = a$ in the $C^*$-norm. Since $\hat{\pi}$ is continuous (as it is a representation), and composition is continuous in $B(H):$ \begin{align*}[\hat{\pi}(a), 2P_N-1] &=  \hat{\pi}(\lim_{n \to \infty}(a_n))(2P_N-1) - (2P_N-1)\hat{\pi}(\lim_{n \to \infty}(a_n))\\ &= \lim_{n \to \infty}(\hat{\pi}(a_n))(2P_N-1) - (2P_N -1)(\lim_{n \to \infty} \hat{\pi}(a_n))\\ &= \lim_{n \to \infty} [\hat{\pi}(a_n), P_N]  \in K(\ell^2(\Z) \otimes H),\end{align*} by the closedness of the compact operators. 
\end{proof}
\begin{prop}\label{CPindex}
    Suppose $A \cong C(X)$ for $X$ compact and Hausdorff and $\vp: X \to X$ a homeomorphism. If $(\ell^2(\Z) \otimes H, \hat{\pi}, 2 P_N -1)$ is a cycle over $C(X) \rtimes_\vp \Z$ as in Proposition \ref{CPodd} and $u \in C(X) \rtimes_\vp \Z$ is the unitary implementing $\vp.$ Then $ \langle [u], [\ell^2(\Z) \otimes H, \hat{\pi}, 2 P_n -1] \rangle = -|N|.$ By choosing to project onto $m <0,$ we obtain classes whose pairing with $[u]$ equals $|N|$ (compare to \cite[Example 2]{H}).
\end{prop}
\begin{proof}
     Observe that: \begin{align*}
          \langle [u], [\ell^2(\Z) \otimes H, \hat{\pi}, 2 P_n -1] \rangle & = \text{F-Index }(P_N \hat{\pi}(u) P_N): P_N(\ell^2(\Z) \otimes \ell^2(Y)) \to P_N(\ell^2(\Z) \otimes \ell^2(Y))\\
          &= \text{F-Index }(e_n \otimes \de_\mu \mapsto e_{n+1} \otimes \de_\mu)\\ &= 0 - \text{dim}(\{e_n \otimes \de_\mu | n =1, \mu \in N\}) \\ &= -|N|.
     \end{align*}  
\end{proof}
We now apply Propositions \ref{CPodd} and \ref{CPindex} to the case of Cantor minimal systems.
\begin{Cor}\label{oddcantor}
    Let $X = \Omega^\N \text{ or } \Omega^\Z$ be the Cantor set for a finite alphabet $\Omega.$ Let $\vp: X \to X$ be a self-homeomorphism and also $\vp:C(X) \to C(X)$ the induced automorphism. Let $N \sub Y$ be a finite set of finite words that occur in $X$ and $\tau:Y \to X$ a function such that $\tau (\mu) \in C_\mu$ for each $\mu \in Y.$ If $\pi_\tau: C(X) \to B (\ell^2(Y))$ is as in Proposition \ref{BPST}, $\hat{\pi}_\tau: C(X) \rtimes_\vp \Z \to B(\ell^2(\Z \times Y))$ is as in Proposition \ref{hat}, and $P_N$ is projection onto the subspace spanned by those $e_n \otimes \de_\mu$ where $n \geq 0, \mu \in N,$ then $(\ell^2(\Z \times Y), \hat{\pi}_\tau, 2P_N-1)$ is an odd Fredholm module that is $p$-summable for all $p >0$ on $C_c(\Z, C(X)).$ If $\vp$ is minimal, the pairing of this cycle with $[u] \in K_1(C(X)\rtimes_\vp \Z)$ equals $ -|N|.$
\end{Cor}
Observe that, by projecting onto the negative instead of positive subspace of $\ell^2(\Z)$ we obtain cycles whose pairing with $u$ is positive, i.e.
\begin{Cor}
    Let $X, \vp, \tau, N$ be as in the previous corollary. Define $P'_N$ to be projection onto $e_n \otimes \de_\mu$ where $n \leq 0, \mu \in N.$ Then $(\ell^2(\Z \times Y), \hat{\pi}_\tau, P'_N)$ is a Fredholm module that is $p$-summable for all $p>0$ on $C_c(\Z,C(X))$ and pairs with $[u] \in K_1(C(X) \rtimes_\vp \Z)$ to $|N|.$ 
\end{Cor}

\begin{Cor}
Let $(H,\rho,F)$ be a cycle as in Proposition \ref{oddcantor}. Suppose $\langle [H,\rho,F], [u] \rangle = z.$ Then the homomorphism $- \otimes_{C(X) \rtimes \Z} [H, \rho,F]: K_0(C(X) \rtimes \Z) \oplus K_1(C(X) \rtimes \Z) \to K_1(\C) \oplus K_0(\C)$ is the map that takes the generator $[u] \in K_1(C(X) \rtimes \Z) \cong \Z$ to \[[Id_z] = z[Id_1] \in K_0(\C) \cong \Z \cong \langle [Id_1] \rangle\] where $Id_z \in M_z(\C)$ is the $z \times z$ identity matrix.
\end{Cor}
\begin{proof}
The index pairing $\langle \cdot, \cdot \rangle: KK_1(\C, C(X) \rtimes \Z) \times KK_1(C(X) \rtimes \Z, \C) \to KK_0(\C,\C)$ coincides with the Kasparov product and, since $K_1(\C) \cong \{0\},$ the non-trivial part of the homomorphism exists only from the odd part of the $K$-theory. 
\end{proof}
Thus, applying Lemma \ref{Cantorodd} to Lemmas \ref{IndexUCT}, \ref{IndexUCTCor} we obtain:
\begin{theorem}\label{oddthm}
    Let $x \in K^1(C(X) \rtimes_\vp \Z).$ Pairings of $x$ with elements of $K_1(C(X) \rtimes_\vp \Z)= \langle [u] \rangle$ can be computed using Connes' trace formula for odd $p>0$ summable cycles.
\end{theorem}

We can also lift the cycles of Proposition \ref{CPodd} to unbounded ones:
\begin{Cor}
The cycles of Proposition \ref{CPodd} can be lifted (along the map $(H, \rho, D) \mapsto (H, \rho, D|D|^{-1})$) to a spectral triple that is $p$-summable for all $p >0$ over $(C_c(\Z,C(X))$ with \[\text{Dom}(D) = \left\{(\xi_n(\mu))| \sum_{n \in \Z, \mu \in Y} W^{|n|+|\mu|} \|\xi_n(\mu)\|^2 < \infty \right\}, \: \: \: D(e_n \otimes \de_\mu) = \begin{cases} W^{|n|+|\mu|}  & n >0, \mu \in N  \\  -W^{|n| + |\mu|} & \text{ else }\end{cases} \] by selecting $W$ large enough (e.g. if $X = \Omega^\N, W > |\Omega|.$)
\end{Cor}
\begin{proof}
We first show that $(\ell^2(\Z) \otimes H,\hat{\pi}, D)$ is a spectral triple. Given that $D$ acts as a diagonal operator with positive real eigenvalues with respect to the basis $\{e_n \otimes \de_\mu\},$ we have that $D$ is symmetric and self-adjoint. Because the commutators with the bounded cycles are finite rank, they remain finite rank when the bounded operator $F$ is lifted to the unbounded operator $D.$ Thus, for $f \in C_c(\Z, C(X)), [f,D]$ is bounded. We claim that $D|D|^{-1} = 2P_N-1.$  To see this, observe that $D$ is diagonal. Thus, \[D|D|^{1}(e_n \otimes \de_\mu) = \begin{cases} 1 & n>0, \mu \in N \\ -1 & \text{ else} \end{cases} = 2 P_N -1.\] Plus, \[\text{Tr}(1+D^2)^{-\f p2} = \sum_m \sum_n (1 + W^{2(|\mu|+|n|)})^{-\f p2} < \infty,\] so long as $p >0$ and $W$ is large enough (e.g. if $X = \Omega^\N, W$ needs to be larger than the size of the alphabet defining $X)$ (as in Theorem \ref{spectripodometer}).
\end{proof}
\begin{Ex}
We return to Example \ref{Bratexamp}. Recall that, in this case, $Y$ is the set of finite paths on the diagram and $X$ is the infinite path space. We can define, for example: $\tau: Y \to X$ via \[\tau(\mu) = \begin{cases} \mu000\ldots & \text{ if }\mu \text{ ends in } 0 \text{ or } 1\\ \mu 1000\ldots & \text{ otherwise }\end{cases}.\] Then for any finite set of finite paths, $N \sub Y,$ we can define a Fredholm module $(\ell^2(\Z \times Y), \hat{\pi}_\tau, 2 P_N-1)$ as in Corollary \ref{oddcantor}. The pairing of this Fredholm module with $[u] \in K_1(C(X) \rtimes_\vp \Z)$ equals $-|N|.$

\end{Ex}
\section{$K$-homology and Spectral Triples on the Cantor Set}
\label{k-hom_spec_trip}
In this section, we exhaust the even $K$-homology of the Cantor set using unbounded cycles (the odd $K$-homology is trivial). We use these cycles to prove uniform summability for odometers in the next section. We show that, if $A$ is an AF-algebra, $K^0(A)$ can be exhausted by finitely summable unbounded cycles. We do so by taking the bounded cycles constructed in \cite[Section 4.2]{R} and lifting them to unbounded ones. We then construct explicit cycles for the Cantor set. We accomplish this by taking the Belissard-Pearson spectral triples of \cite{BP} and modifying them to exhaust $K^*(C(X)).$ (Note that $K^1(C(X)) \cong \{0\}$). We then show that these cycles agree, under the bounded transform, with the even cycles for an AF-algebra constructed in \cite[Section 4.2]{R}. 

\subsection{Unbounded Cycles on AF-algebras}
\begin{prop}{(Unbounded Version of \cite[Theorem 4.2.1]{R})}\label{UnboundedAF}
    Let $A = \br{\bigcup_n A_n}$ be an AF-algebra where $\{A_n\}_{n=1}^\infty$ is a decomposition such that each $A_n$ is finite-dimensional. Then for every class $x \in K^0(A)$ and $p > 0,$ there is an unbounded Fredholm module that is $p$-summable on $\bigcup_n A_n$ that represents $x.$ The bounded transform of such a module is $q$-summable for all $q>0$ summable on $\bigcup_n A_n$ (independent of choice of $p$).
\end{prop}
\begin{proof}
 Suppose $I \in \Hom(K_0(A), \Z).$ We construct an unbounded module whose class, $x \in K^0(A)$ satisfies $\langle y, x \rangle = I(x)(y) \in \Z$ for all $y \in K_0(A).$ We utilize the construction in \cite[Theorem 4.2.1]{R} to find a Hilbert space and representation for our unbounded cycle. We construct an unbounded operator to produce an unbounded cycle representing the class $x.$ Following \cite[Theorem 4.2.1]{R}, select a separable infinite-dimensional Hilbert space $H.$ Then, as done in \cite[Theorem 4.2.1]{R}, we define $\vp_n^\pm: A_n \to B(H)$ such that each $\vp_{n+1}^\pm$ extends $\vp_n^\pm.$ Further, these $\vp_n^\pm$ satisfy that there is a set of matrix units $\{e_{ij}^{(k)}\}$ for each $A_n$ and pairwise orthogonal projections $P_n^{(k)}, P_n^{\pm (k)} \in B(H)$ for which following properties hold for each $i,k$:
 \begin{enumerate}
\item $\vp_n^\pm (e_{11}^{(k)}) = P_n^{(k)} + P_n^{\pm (k)}$ 
\item $\vp_n^+(e_{i1}^{(k)})P_n^{(k)} = \vp_n^-(e_{i1}^{(k)})P_n^{(k)}$
\item $\text{rank } P_n^{(k)} = \infty, \text{rank } (P_n^{\pm^{(k)}}) < \infty $
\item $\begin{cases} \text{rank } P_n^{+(k)} = I(e_{11}^{(k)}), P_n^{-(k)} = 0 & \text{ if } I(e_{11}^{(k)}) \geq 0 \\ \text{rank } P_n^{-(k)} = |I(e_{11}^{(k)})|, P_n^{+(k)} = 0 & \text{ if } I(e_{11}^{(k)}) < 0 \end{cases}$.
 \end{enumerate}
    Having done so, we can define the even unbounded Fredholm module $(H \oplus H, \varinjlim \vp_n^+ \oplus \varinjlim \vp_n^-, D'),$ which differs from Rave's construction only in the operator. We define the operator $D: \bigcup_n \bigcup_k \text{range } P_n^{(k)} \cup P_n^{\pm (k)} \to \text{range } P_n^{(k)} \cup P_n^{\pm (k)}.$ Note that $\varinjlim \vp_n^+ \oplus \varinjlim \vp_n^- (A_n) = \bigcup_k \text{ range } P_n^{(k)} \cup \text{ range } P_n^{+(k)} \cup \text{ range } P_n^{- (k)}.$ Thus, by \cite[Theorem 2.1 (ii)]{CI}, for each $r > 0$ there is a sequence $\alp_n$ such that, letting \[D = \sum_{n=1}^\infty \alp_n \left(\sum_k P_n^{(k)}+P_n^{+ (k)} + P_n^{-(k)} - \sum_{l=1}^{n-1} \sum_k P_l^{(k)} + P_{l}^{+ (k)} + P_{l}^{- (k)}\right)\] and \[D' = \begin{bmatrix} 0 & D\\ D & 0 \end{bmatrix},\] gives $D'$ not only compact resolvent but also $L^p(H)$ resolvent for $p \geq r.$ Next, suppose $a \in \bigcup_{n=1}^\infty A_n.$ Then $(\varinjlim \vp_n^+ \oplus \varinjlim \vp_n^- (a)) \bigcup_n \bigcup_k \text{ range } P_n^{(k)} \cup P_n^{\pm (k)} \sub \bigcup_n \bigcup_k \text{ range } P_n^{(k)} \cup P_n^{\pm (k)},$ as $\varinjlim \vp_n^+ \oplus \varinjlim \vp_n^- (A_n) = \bigcup_k \text{ range } P_n^{(k)} \cup \text{ range } P_n^{+(k)} \cup \text{ range } P_n^{- (k)}.$ Thus, the representation preserves the domain of $D'.$ Additionally, the commutators of $D'$ with elements of $\bigcup_n A_n,$ by construction, are finite rank. Further, under the bounded transform $(H, \rho, D) \mapsto (H,\rho, D|D|^{-1})$ this cycle maps to $\left(H \oplus H, \varinjlim \vp_n^+ \oplus \varinjlim \vp_n^-, \begin{bmatrix} 0 & 1 \\ 1 & 0 \end{bmatrix} \right),$ which agrees with Rave's construction since we have constructed the Hilbert space and representation using his method. Per Rave, the index map associated to the bounded transform of this cycle is $I.$ This bounded cycle is $q$-summable for all $q >0$ since the commutators are finite rank and $\begin{bmatrix} 0 & 1 \\ 1 & 0 \end{bmatrix}$ is a self-adjoint unitary. Recall that, \cite[Section 4.1]{R}, the index maps on $K_0(A)$ determine the $K$-homology $K^0(A).$ Thus, for any $p > 0$ we can attain any index map on $K_0(A)$ with a spectral triple that is $p$-summable on $\bigcup_n A_n,$ and thus every class in $K^0(A)$ as such.
\end{proof}

\subsection{Rave's bounded cycles for the Cantor Set}
We seek to use Proposition \ref{UnboundedAF} to write explicit unbounded cycles to exhaust $K^0(C(X))$ for $X$ the Cantor Set. To do so, we write down bounded cycles for $K^0(C(X))$ from \cite[Section 4.2]{R}  and show how they lift to unbounded ones. To begin, we describe the well-known AF-structure of $C(X)$ when $X$ is a Cantor set.
\begin{Def}
    Let $\chi_Z$ denote the \textbf{indicator function} on the set $Z.$ Observe that, when $Z$ is a cylinder set (i.e.  $Z = C_\mu,$ for some $\mu \in Y), \: \chi_{C_\mu}$ is locally constant and thus is continuous.
\end{Def}  
\begin{prop}{\cite{GM}}\label{AFfiltration}
    An AF filtration of $C(X) = A$ is given by \[A_n = \C^{|\{\mu |\: \mu \in Y, \:|\mu| = n\}|} \cong \oplus_{|\mu| = n} \C [\chi_{C_\mu}]\] with the inclusion maps $\iota: A_n \to A_{n+1}$ given by the partition $C_\mu = \bigcup_{\lambda \in Y, |\lambda| =1, \mu \lambda \in Y} C_{\mu \lambda}.$ Thus $C(X) = \br{\bigcup_n A_n}.$
\end{prop}

\begin{prop}{\cite[Section 4.2]{R}} 
    The map $\gm:K^0(C(X)) \to \text{Hom}(K_0(C(X)), \Z)$ given by $\gm(x)(-):= \langle -,x \rangle$ is an isomorphism.
\end{prop}
 \begin{prop}{\cite[Theorem 4.2.1]{R}}\label{Rave}
 Let $x \in K^0(C(X)).$ Then there is an even Fredholm module that is $p$-summable for all $p>0$ on $\text{span} \{\chi_{C_\mu}|\mu \in Y\}$ that represents $x.$ 
 \end{prop} 
This is the application of \cite[Theorem 4.2.1]{R} to $C(X)$ with AF-filtration given by Proposition \ref{AFfiltration}.
\subsection{Belissard-Pearson Spectral Triples}
We now lift the cycles of Proposition \ref{Rave} to unbounded ones. We know we can lift to summable cycles from Proposition \ref{UnboundedAF}. We start with the even Belissard-Pearson Spectral triples associated to weak choice functions, as detailed in \cite{GM} and \cite{BP}. Towards this end, 
\begin{Def}
    Let $\Omega$ be a finite set (of symbols), which we call an \textbf{alphabet}.
\end{Def} 
Let the Cantor Set $X$ be modeled by $\Omega^\N$ as in Example \ref{Cantoralph}.

\begin{Def}
    Let $Y$ be the set of finite words that appear in $\Omega^\N,$ which we call the \textbf{language} of $X$. Observe that $Y$ is countable.
\end{Def} 

\begin{Def}
    For $\mu \in Y,$ let $ C_\mu$ denote the \textbf{cylinder set} of sequences in $\Omega^\N$ that begin with the word $\mu,$ i.e. $C_\mu =\{(a_1,a_2,\ldots) \sub \Omega^{\N}| \: a_1a_2\ldots a_{|\mu|} = \mu\}.$ Endow $\Omega^\N$ with the topology that has basis given by these cylinder sets. 
\end{Def}
\begin{Def}\label{choice} Let $\tau = (\tau_+,\tau_-): Y \to X \times X.$ We say $\tau$ satisfies the \textbf{cylinder condition} if, for each $\mu \in Y, \: \tau_+(\mu), \tau_-(\mu) \in C_\mu.$ If $\tau$ satisfies the cylinder condition, we call $\tau$ a \textbf{weak choice function}.
\end{Def}
\begin{prop}{\cite{BP}}\label{pitau}
    For each weak choice function, $\tau,$  we can define a representation \[\pi_{\tau} = \pi_{\tau_+} \oplus \pi_{\tau_-}: C(X) \to B(\ell^2(Y) \oplus \ell^2(Y))\] via: \[\pi_{\tau_\pm}(f)(\xi) (\nu) = f(\tau_\pm(\nu)) \cdot \xi(\nu)\] for $\nu \in Y, \xi \in \ell^2(Y), f \in C(X).$ 
\end{prop}
Next, define the odd operator $D: C_c(Y) \oplus C_c(Y) \to \ell^2(Y) \oplus \ell^2(Y)$ via \[D \left(\begin{bmatrix} \xi_0 \\ \xi_1 \end{bmatrix}\right)(\nu) := e^{|\nu|} \begin{bmatrix} \xi_1(\nu) \\ \xi_0(\nu)\end{bmatrix} .\] From \cite{BP},\cite{GM} we have that:
\begin{prop}\label{BPST}
    For each weak choice function, $\tau, (\ell^2(Y) \oplus \ell^2(Y), \pi_\tau, D)$ forms a spectral triple over $C(X).$ The Lipschitz algebra of any such spectral triple is the algebra of Lipschitz functions on $X.$ 
\end{prop}
\begin{prop}
    Let $\tau$ be a weak choice function. Let $p >0.$ Then there is an operator $D_p: C_c(Y) \oplus C_c(Y) \to \ell^2(Y) \oplus \ell^2(Y)$ given by \[D_p \left(\begin{bmatrix} \xi_0 \\ \xi_1 \end{bmatrix}\right) (\nu) := A_p(|\nu|) \begin{bmatrix} \xi_1(\nu) \\ \xi_0(\nu) \end{bmatrix}\] such that $(\ell^2(Y) \oplus \ell^2(Y), \pi_\tau, D_p)$ is a spectral triple that is $p$-summable on the algebra of Lipschitz functions (\cite[Remark 4.1.5]{GM}).
\end{prop}

\begin{lemma} Consider the bounded transform for spectral triples with invertible operator $D, \: (H, \rho, D) \mapsto (H, \rho, D|D|^{-1}).$ Per \cite{GM}, the Belissard-Pearson spectral triples map to even cycles of Proposition \ref{Rave} under this map.
\end{lemma}
\begin{proof} First, suppose $\tau = (\tau_+,\tau_-)$ is a weak choice function and $(\pi_\tau \oplus \pi_\tau, \ell^2(Y) \oplus \ell^2(Y), D)$ is a Belissard-Pearson Spectral Triple as in Proposition \ref{BPST}. Observe that the bounded transform of the Belissard-Pearson Dirac Operator is the operator $F = \begin{bmatrix} 0 & 1 \\ 1 & 0 \end{bmatrix}.$ Further, the representation $\pi_\tau$ is built up inductively on levels of the AF-filtration. Specifically, using the notation of Proposition \ref{Rave}, \[\vp_n^\pm(\chi_{C_\mu}) = \pi_{\tau_\pm}(\chi_{C_\mu}) = P_{\tau_+^{-1}(C_\mu) \cap \tau_-^{-1}(C_\mu)} + P_{\tau_\pm^{-1}(C_\mu) \setminus \tau_\mp^{-1}(C_\mu)},\] where $P_Z$ denotes the projection onto the set $Z.$ To see that properties (1-4) of Proposition \ref{Rave} are satisfied, notice that \[\pi_{\tau_+}(\chi_{C_\mu}) P_{\tau_+^{-1}(C_\mu) \cap \tau_-^{-1}(C_\mu)} = \pi_{\tau_-}(\chi_{C_\mu}) P_{\tau_+^{-1}(C_\mu) \cap \tau_-^{-1}(C_\mu)}  = P_{\tau_+^{-1}(C_\mu) \cap \tau_-^{-1}(C_\mu)}.\] Also, observe that $P_{\tau_+^{-1}(C_\mu) \cap \tau_-^{-1}(C_\mu)}$ is infinite rank. This is because $\nu \in \tau_+^{-1}(C_\mu) \cap \tau_-^{-1}(C_\mu)$ for each $\nu \in Y$ that begins with $\mu.$ Since $\mu$ is in at least one infinite word, the set of such $\nu$ is infinite. Then, for each $\mu, \: P_{\tau_\pm^{-1}(C_\mu) \setminus \tau_\mp^{-1}(C_\mu)}$ is finite rank, as $\tau_\pm^{-1}(C_\mu) \setminus \tau_\mp^{-1}(C_\mu)$ only contain subwords of $\mu.$ Finally, for property (4), observe that, for $\mu \in Y,$ \begin{align*}\gm(x)([\chi_{C_\mu}]) &= \text{F-Index}(\pi_{\tau_-}(\chi_{C_\mu})\pi_{\tau_+}(\chi_{C_\mu})):\pi_{\tau_+}(\chi_{C_\mu})\ell^2(Y) \mapsto \pi_{\tau_-}(\chi_{C_\mu})\ell^2(Y)) \\ &= \text{Rank } P_{\tau_+^{-1}(C_\mu) \setminus \tau_-^{-1}(C_\mu) } - \text{Rank } P_{\tau_-^{-1}(C_\mu) \setminus \tau_+^{-1}(C_\mu)} .\end{align*}
The generators of $K_0(C(X))$ at the $n$-th level of the AF-filtration are the indicator functions on cylinder sets defined by words of length $n.$ For each of these indicator functions, the positive and negative parts of the representation agree up to a finite rank projection. When paired with these indicator functions, the index is equal to the rank of the difference of these projections. This is computed in both \cite{R} and \cite{GM}.
\end{proof}
\subsection{Restricted Cycles and the rest of the $K$-homology}
Recall that, for $X$ the Cantor set, the map $\gm: K^0(C(X)) \to \text{Hom}(K_0(C(X)), \Z)$ is an isomorphism. We use this result to exhaust $K^0(C(X)).$ To do so, we take an arbitrary index map and construct an unbounded cycle whose class gives said map. 
Observe that, 
\begin{prop}{(proof in \cite[Theorem 4.2.1]{R})}
    The set $\text{span}\{[\chi_\mu]| \: \mu \in Y\}$ is dense in $K^0(C(X)) \cong C(X,\Z).$ Since index maps are continuous, they are determined on this set.
\end{prop}
To attain all index maps, we must allow for modifications of the Belissard-Pearson cycles. There are two restrictions on these cycles we need to account for. First,
\begin{lemma}\label{obst1}
The $K$-homology class of a Belissard-Pearson spectral triple from a function that satisfies the cylinder condition (see Definition \ref{choice}) pairs with the class of $[1_X] \in K_0(C(X))$ to zero.
\end{lemma} 
\begin{proof}
This follows from the computation in the proof of \cite[Lemma 4.3.4]{GM}. Let $\tau = (\tau_+, \tau_-): Y \to X \times X$ be a function that satisfies $\tau_+(\mu), \tau_-(\mu) \in C_\mu$ for each $\mu \in Y.$ Let $x = (H,\pi_\tau,D)$ be a Belissard-Pearson spectral triple associated to $\tau$ (see Proposition \ref{BPST}). Let $\gm(x) : K_0(C(X)) \to \Z$ be given by $\langle -, [\pi_\tau, H,D] \rangle = \langle -, [\pi_\tau, H, D|D|^{-1}] \rangle$ (see Proposition \ref{BT}). For any word $\nu \in Y,$ the function $f \equiv 1$ takes the value $1$ on both $\tau_+(\nu)$ and $\tau_-(\nu).$ This means that, for any $n\geq 1, \ds \sum_{|\mu| = n} \gm(x)([\chi_{C_\mu}]) =0.$
\end{proof}
In addition, there is an upper bound on values the index map can take based on the length of the cylinder sets:
\begin{lemma} \label{obst2}
Suppose $\tau$ is a weak choice function that satisfies the cylinder condition (see Definition \ref{choice}). Let $(H,\pi_\tau,D)$ be a Belissard-Pearson spectral triple associated to $\tau.$ Let $\mu \in Y.$ Then $|\langle [\chi_{C_\mu}], [H,\pi_\tau,  D|D|^{-1}] \rangle | \leq |\mu|.$
\end{lemma}
\begin{proof} This follows from the computations in \cite[Section 4]{GM}. Specifically, for $\mu \in Y,$ denote by $S_\mu = \{ \nu | \mu = \nu \lambda \text{ for some } \lambda \in Y\}.$ Then $|S_\mu| = |\mu|.$ Further, \[\langle [\chi_{C_\mu}], [x] \rangle = \#\{\nu \in S_\mu | \tau_+(\nu) \in C_\mu, \tau_-(\nu) \not \in C_\mu\} - \#\{\nu \in S_\mu| \tau_-(\nu) \in C_\mu, \tau_+(\nu) \not \in C_\mu\}.\]
\end{proof}
However, by altering the weak choice functions and taking direct sums, we can construct cycles to obtain any index map in $\Hom(K_0(C(X)), \Z).$ 

To obtain cycles whose classes pair with the identity non-trivially, we alter the representations allowed. Specifically, we start by selecting a word $\mu.$ Then we alter the representation $\pi: C(X) \to B(\ell^2(Y))$ by restricting to the set $\chi_\mu$ before multiplying. Recall that \cite{GM}, \cite{BP} define the representation $\pi_\tau: C(X) \to B(\ell^2(Y))$ as the composition of the pullback $\tau^*: C(X) \to C_b(Y)$ with the multiplication representation $\rho: C_b(Y) \to B(\ell^2(Y)).$ Then:
\begin{prop}\label{restrep}
    For each $\mu \in Y,$ and $\tau = (\tau_+, \tau_-)$ a weak choice function, we obtain a representation $\pi_{\tau_\mu} := \rho \circ \chi_\mu \circ \tau^*$ where the map $\chi_\mu: C(X) \to C(X)$ takes $f$ to $f\cdot \chi_\mu.$ 
\end{prop}
\begin{proof}
    This is a representation as $\pi_{\tau_\mu}(f) = \pi_\tau(f \cdot \chi_\mu)$ and $\chi_\mu$ is a projection.  
\end{proof}
\begin{prop}
    Suppose $\tau = (\tau_+, \tau_-)$ is a weak choice function, $\pi_{\tau_\mu}$ is as in Proposition \ref{restrep} and $D$ is as in Proposition \ref{BPST}. Then the computations in \cite{BP} show that, for a Lipschitz function $f, [\pi_{\tau_\mu}(f), D]$ is bounded. 
\end{prop}
\begin{proof}
    Let $\mu \in Y.$ Then \[\# S_\mu = \{\nu \in Y|\tau_+(\nu) \in C _\mu , \tau_-(\nu) \not \in C_\mu \text{ or } \tau_-(\nu) \in C_\mu, \tau_+(\mu) \not \in C_\mu\}\] is finite. For $\nu \not \in S_\mu [\pi_{\tau_\mu}(f), D](\nu) = [\pi_{\tau}(f),D](\nu).$ By, \cite{BP}, $[\pi_\tau(f),D]$ is bounded if $f$ is Lipschitz. Thus $[\pi_{\tau_\mu}(f), D] = [\pi_{\tau}(f),D)] + F$ where $F$ is finite rank, so $[\pi_{\tau_\mu}(f), D]$ is bounded.
\end{proof}
\begin{Cor}\label{RBPST}
For a weak choice function $\tau$ and $\mu \in Y,$ $(\ell^2(Y) \oplus \ell^2(Y), \pi_{\tau_\mu}, D)$ is an unbounded Fredholm module on $\text{Lip}(X).$ We will refer to these modules as \textbf{Restricted Belissard-Pearson Cycles}.
\end{Cor}
\begin{proof}
    Observe that we have only altered the representation. Thus, it remains to check that the new representation is indeed a representation and that commutators of the images of Lipschitz functions under this representation with the operator $D$ are bounded. We have done so in the previous two propositions.
\end{proof}
\textit{Note:} Belissard-Pearson spectral triples are restricted Belissard-Pearson cycles with $\mu = \ee.$
\begin{prop}
If $(\ell^2(Y) \oplus \ell^2(Y),\pi_\tau,D)$ is a Belissard-Pearson spectral triple for a weak choice function $\tau$ and $\mu \in Y,$ then \[\langle [1_X], [\ell^2(Y) \oplus \ell^2(Y),\pi_\tau,D|D|^{-1}] \rangle = \langle [\chi_{C_\mu}], [\ell^2(Y) \oplus \ell^2(Y),\pi_{\tau_\mu},D|D|^{-1}] \rangle.\]
\end{prop} 
\begin{proof}
    In particular, \[\langle [1_X], [\ell^2(Y) \oplus \ell^2(Y),\pi_{\tau},D|D|^{-1}] \rangle = \langle [\chi_{C_\mu}], [\ell^2(Y) \oplus \ell^2(Y),\pi_{\tau_\mu},D|D|^{-1}] \rangle. \] This is because: \begin{align*}
        \langle [\chi_{C_\mu}], [\ell^2(Y) \oplus \ell^2(Y),\pi_{\tau},D|D|^{-1}] \rangle &= \text{F-Index}(\pi_{\tau_+}(\chi_{C_\mu})\pi_{\tau_-}(\chi_{C_\mu})\\ &=  \text{Rank } P_{\tau_+^{-1}(C_\mu) \setminus \tau_-^{-1}(C_\mu) } - \text{Rank } P_{\tau_-^{-1}(C_\mu) \setminus \tau_+^{-1}(C_\mu)}\\ &=\text{Rank } P_{\tau_{\mu_+}^{-1}(X) \setminus \tau_{\mu_-}^{-1}(X) } - \text{Rank } P_{\tau_{\mu_-}^{-1}(X) \setminus \tau_{\mu_+}^{-1}(X)}\\ &= \langle [1_X],[\ell^2(Y) \oplus \ell^2(Y),\pi_{\tau_\mu},D|D|^{-1}] \rangle.\end{align*}
\end{proof}
\begin{prop}
Let $A$ be a finite alphabet. Let $X \sub A^\N$ be a Cantor Set. Let $Y$ be the set of finite words in $X.$ Suppose $x \in K^0(C(X))$ is such that $\gm(x)([\chi_{C_\mu}]) := \langle [\chi_{C_\mu}], x \rangle$ satisfies $|\gm(x)([\chi_{C_\mu}])| < |\mu|$ for each $\mu \in Y$ and $[\chi_{C_\mu}] \in K_0(C(X)).$ Then there exists a Restricted Belissard Pearson cycle (see Proposition \ref{RBPST}) $(\ell^2(Y,\C^2), \pi_{\tau_\mu}, D)$ ($\mu$ could be $\ee$) representing $x.$
\end{prop}
\begin{proof}
Recall that $\gm(x)$ is linear. Observe that, for each $\mu,$ there exist words $\nu_1, \ldots, \nu_k$ of length $|\mu|+1$ such that $\chi_{C_\mu} = \sum_{i=1}^k \chi_{C_{\nu_i}}.$ Thus $\gm(x)([\chi_{C_\mu}]) = \sum_{i=1}^k \gm(x)([\chi_{C_{\nu_i}}].$ Since $\gm(x)$ is determined on these cylinder sets, we can determine $\gm(x)$ by the sequence of integers \[\gm(x)([1_X]), \: \gm(x)([\chi_{C_0}]), \: \gm(x)([\chi_{C_{00}}]), \: \gm(x)([\chi_{C_{10}}]), \ldots\] obtained by evaluating $\gm(x)$ on words that end in $0.$ To obtain a class that gives such a map, begin with any function $\tau_-: Y \to X$ that satisfies the cylinder condition of \cite{GM}. Now, select a word $\mu,$ such that $|\mu| = |\gm(x)([1_X])| +1$ that maximizes $|\gm(x)([1_{\chi_{C_\nu}}])|$ amongst words $\nu$ such that $|\nu| = |\gm(x)([1_X])| + 1.$ We define $\tau = (\tau_{+}, \tau_{-})$ and the representation $\pi_{\tau_\mu}$ by defining for each $\nu \in Y,$ $\tau_+(\nu)$ so that, for each $\mu \in Y,$ \[\#\{\nu| \tau_+(\nu) \in C_\mu, \tau_-(\nu) \not \in C_\mu\} -  \#\{\nu| \tau_-(\nu) \in C_\mu, \tau_+(\nu) \not \in C_\mu\} = I ([\chi_{C_\mu}]).\] We can do so, so long as, for each $\mu \in Y, |\gm(x)(\chi_{C_\mu})| < |\mu|.$ Note that a cylinder set determined by a word of length $n$ can be partitioned into cylinder sets of longer length. Thus, the index of classes of indicator functions on cylinder sets of length $n$ must be the sum of the indicies of classes of indicator functions on the sets in such a partition. This relation is precisely the relation in $K_0(C(X)) \cong C(X,\Z).$ i.e. if $\chi_{C_{\nu_1}} + \chi_{C_{\nu_2}} = \chi_{C_\mu}$ then $\gm(x)([\chi_{C_{\nu_1}}]) +\gm(x)([\chi_{C_{\nu_2}}]) = \gm(x)([\chi_{C_\mu}])$ because \[ \#\{\nu| \tau_+(\nu) \in C_\mu, \tau_-(\nu) \not \in C_\mu\} =  \#\{\nu| \tau_+(\nu) \in C_{\nu_1}, \tau_-(\nu) \not \in C_{\nu_1}\} + \#\{\nu| \tau_+(\nu) \in C_{\nu_2}, \tau_-(\nu) \not \in C_{\nu_2}\}.\] The same holds when we reverse the roles of $\tau_+$ and $\tau_-.$
\end{proof}
Now, we use these results and Proposition \ref{UnboundedAF} to exhaust $K^0(C(X))$ with unbounded cycles. For any class $x \in K^0(C(X))$ we construct an unbounded cycle obtained as a (perhaps infinite) `linear combination' of these restricted cycles with the Belissard-Pearson spectral triples to represent $x.$ We do so to overcome the obstacles of Propositions \ref{obst1}, \ref{obst2}. We know the cycle constructed represents $x$ because we check that its image under $\gm$ coincides with $\gm(x) \in \text{Hom}(K_0(C(X)), \Z).$ In order to account for those $x \in K^0(C(X))$ with infinitely many $\mu$ such that $|\gm(x)([\chi_{C_\mu}])| \geq |\mu|,$ we take the Hilbert Space to be $\ell^2(\N) \otimes \ell^2(Y,\C^2).$ Towards this end, suppose $X \sub \Omega^\N$ for a finite alphabet $\Omega.$ 
\begin{Cor}\label{UnboundedCantor}
    Let $\Omega$ be a finite alphabet, $X \sub \Omega^\N$ a Cantor set,  and $[x] \in K^0(C(X))$. Let $p > 0$ and define $W$ so that $W^p > |\Omega|.$ Then we can select a word $\nu \in Y,$ a weak choice function $\tau,$ a sequence of weak choice functions $(\tau_n)_{n \in \N},$ and Hilbert Space $H = \ell^2(\Z_{\geq 0}) \otimes \ell^2(Y,\C^2).$ Then, define the representation $\tilde{\pi} = \tilde{\pi_+} \oplus \tilde{\pi}_-: C(X) \to B(H)$ by \begin{align*} \tilde{\pi}(f)\left(\de_n \otimes \begin{bmatrix} \xi \\ \xi' \end{bmatrix}\right) &= \left(\de_n \otimes \begin{bmatrix} f(\tau_{n_+}(\xi)) \xi \\ f(\tau_{n_-}(\xi')) \xi' \end{bmatrix} \right) \text{ if } n \geq 1 \text{ and } \\ \tilde{\pi}(f)\left(\de_0 \otimes \begin{bmatrix} \xi \\ \xi' \end{bmatrix}\right) &= \left(\de_0 \otimes \begin{bmatrix} f(\tau_{\mu_+}(\xi)) \xi \\ f(\tau_{\mu_-}(\xi')) \xi' \end{bmatrix} \right).\end{align*} We obtain a $p$-summable unbounded Fredholm module $(H, \tilde{\pi},D) $ with $D: \text{Dom}(D) \to H$ given by \[D\left(\de_n \otimes\begin{bmatrix}  \xi \\ \xi' \end{bmatrix} (\nu)\right) = W^n \de_n \otimes W^{|\nu|}\begin{bmatrix} \xi'(\nu) \\ \xi(\nu) \end{bmatrix} \text{ over  }\A = \text{span}\{\chi_{C_\mu} | \mu \in Y \}.\] We can make these selections so that the module represents $[x] \in K^0(C(X)).$
\end{Cor}
\begin{proof}
To begin, if $\langle [x], [1_X] \rangle = 0,$ we can choose $\nu = \ee.$ Otherwise, if $\langle [x], [1_X] \rangle = k,$ we choose $\nu$ such that $|\nu| = k +1$ and $\tau_{\nu}$ such that \begin{enumerate}
    \item \[\text{rank } \pi_{\tau_{\nu_+}}(\chi_{C_\mu}) - \text{rank }\pi_{\tau_{\nu_-}}(\chi_{C_\mu}) = \langle [x], [\chi_{C_\mu}] \rangle\] for all $\mu$ such that $|\langle [x], [\chi_{C_\mu}]\rangle| < |\mu|,$ 
    \item for words $\mu$ such that $\langle [x], [\chi_{C_\mu}] \rangle \geq |\mu|, \text{ rank } \pi_{\tau_{\nu_+}}(\chi_{C_\mu}) - \text{ rank }\pi_{\tau_{\nu_-}}(\chi_{C_\mu})= |\mu|-1,$ and 
    \item if $\langle [x], [\chi_{C_\mu}] \rangle \leq -|\mu|, \text{rank } \pi_{\tau_{\nu_+}}(\chi_{C_\mu}) - \text{rank }\pi_{\tau_{\nu_-}}(\chi_{C_\mu})= 1-|\mu|.$
\end{enumerate} 
Then, denote by $\{\mu_i\}$ those $\lambda \in Y$ such that $|\gm(x)([\chi_{C_\lambda}])| > |\lambda|$ ordered by increasing length. First, suppose $\langle [x], [\chi_{\mu_1}] \rangle  = m.$  Let $r > |\f{m}{|\mu_1|}|.$  Define $\tau_1, \ldots, \tau_r$ so that 
\begin{enumerate}
    \item For all $1 \leq j \leq r, \: \tau_{j_-}(\mu) = \tau_{j_+}(\mu) $ for all $\mu$ such that $|\langle [x], [\chi_{C_\mu}] \rangle| < |\mu|$ and  
    \item   \[\text{rank } (\pi_{\tau_{\nu_0^+}} \oplus_{j \leq r} \pi_{\tau_{j^+}} (\chi_{C_{\mu_1}})) - \text{rank } (\pi_{\tau_{\nu_0^-}} \oplus_{j \leq r} \pi_{\tau_{j_-}} (\chi_{C_{\mu_1}})) = \langle x, [\chi_{C_\mu}] \rangle.\]
\end{enumerate}  Then, suppose $\tau_j$ is defined for all $j < l.$ Inductively for each $i,$ we define $\tau_{l}, \ldots \tau_{l+r'}$ so that 
\begin{enumerate}
    \item  For all $l \leq j \leq j+r', \tau_{j_-}(\mu) = \tau_{j_+}(\mu) $ for all $\mu$ such that $|\langle [x], [\chi_{C_\mu}] \rangle| < |\mu|$ or so that $\mu = \mu_s$ for $s < i$ and
    \item $\tau_{l_\pm}, \ldots \tau_{l +r'_\pm}$  satisfy  \[\text{rank } (\pi_{\tau_{\nu_0^+}} \oplus_{j \leq l+r'} \pi_{\tau_{j^+}}(\chi_{C_{\mu_i}})) - \text{rank } (\pi_{\tau_{\nu_0^-}} \oplus_{j \leq l+r'} \pi_{\tau_{j_-}}(\chi_{C_{\mu_i}})) = \langle x, [\chi_{C_\mu}] \rangle.\] 
\end{enumerate} 
In this way, we inductively build up a representation to create a cycle that has the specified index map. Note that, \[\text{Trace}(|D|^{-\f{p}{2}}) = \sum_{k \in \Z_{\geq 0},\: \mu \in Y} W^{-p(k+|\mu|)}  = (2-\frac{1}{W^p}) \sum_{n\in \Z_{\geq 0}} \frac{|\{\mu \in Y | \:  |\mu| = n\}|}{W^{-pn}} \leq  (2-\frac{1}{W^p}) \sum_{n \in \Z_{\geq 0}} \left(\frac{|\Omega|}{W^p}\right)^{n} < \infty.\] Thus, this unbounded Fredholm module is $p$-summable. Additionally, for each $\mu$ and large enough $j, \: \tau_{j_-}(\mu)=\tau_{j_+}(\mu).$ Thus, for each $\mu, [\tilde{\pi}_{\tau}(\chi_{C_\mu}),D] $ is finite rank. This ensures that the bounded transform of this triple will be $d$-summable for all $d>0.$ 
\end{proof}

\section{Odometers} 
Odometers are a particularly nice class of Cantor minimal systems. The crossed-product algebras associated to them are Bunce-Deddens algebras \cite{BuD}. Their $K$-theory is well-understood \cite{P}, \cite{KMP}. From the UCT or Pimsner-Voiculescu sequence, we can compute their $K$-homology. These $C^*$-algebras are not Poincar\'e Duality algebras, and thus the approaches to proving uniform finite summability of \cite{EN}, \cite{GM}, \cite{G} will not work for these algebras. Instead, we use the cycles of the previous section on the Cantor set and extend them to the crossed product algebra for an odometer. We can do so because odometers are \textit{metrically equicontinuous}. Thus, we can use the results of \cite{HSWZ} to obtain unbounded cycles on $C(X) \rtimes_\vp \Z.$ These unbounded cycles exhaust the $K$-homology of $C(X) \rtimes_\vp \Z.$ Specifically, we prove the following, which is the main result of the section.
\begin{prop}
Let $(X,\vp)$ be an odometer. Let \[\A = \text{span}\{\chi_{C_\mu}| \: \mu \in Y \} \text{ and } C_c(\Z,\A):= \left\{\sum_{k=L}^K a_ku^k \: \: L \leq K \in \Z, a_k \in \A \right\} \text{ where } u \text{ implements } \vp.\] Then, for every class $x \in K^*(C(X) \rtimes_\vp \Z)$ and $p >1,$ there is an unbounded Fredholm module that $p$-summable on $C_c(\Z,\A)$ that represents $x.$
\end{prop}
\subsection{Properties of Odometers}
\begin{Def}
    A dynamical system, $(X, \vp),$ is \textbf{equicontinuous} if, for all $\ve > 0,$ there is a $\de >0$ such that, for all $x,y \in X, \: d(x,y) < \de$ implies $d(\vp^n(x), \vp^n(y)) < \ve$ for all $n \in \Z.$
\end{Def}
\begin{prop}{\cite{F} \cite[Proposition 3.9]{HSWZ}}\label{equicontinuosodometers}
    Odometers are precisely, up to topological conjugacy, the Cantor minimal systems that are equicontinuous.
\end{prop}
Odometers additionally are \textit{metrically equicontinuous}, in that:
\begin{prop}{\cite[Proposition 1, Section 2.3]{HSWZ}}
    An equicontinuous action $\vp$ on a metric space $X$ is \textbf{metrically equicontinuous} if there is an equivalent metric such that $\vp$ is isometric.
\end{prop}
This can be seen through:
\begin{prop}{(See \cite{F}, \cite{HSWZ})}
    In the metric of Example \ref{Cantoralph}, odometers are isometric, i.e. $d(x,y) = d(\vp(x),\vp(y))$ for all $x,y \in X.$
\end{prop}

\subsection{$K$-homology for Odometers}
The $K$-theory of odometers has been well-studied. From \cite{P}, for example, we have that, for an odometer $(X,\vp)$ associated to $\{d_i\}, \\ \:K_0(C(X) \rtimes_\vp \Z) \cong \left \{\frac{k}{n_1n_2\cdots n_m}| k, n_1, \ldots, n_m \in \Z, \: m \geq 1 \right\} \sub \Q$ and $K_1(C(X) \rtimes_\vp \Z) \cong \Z.$ Since the $K$-theory of odometers is understood \cite{P}, \cite{KMP}, their $K$-homology can be computed using the UCT. This has been done in \cite{KMP}.
    \begin{prop}{\cite{KMP}}
 If $\vp:X \to X$ is an odometer, then $K^0(C(X) \rtimes_\vp \Z) \cong 0.$ 
    \end{prop}

\begin{prop}{\cite{KMP}}
For an odometer associated to the supernatural number $p_1p_2\cdots $ for primes $p_i,$ \[K^1(C(X) \rtimes_\vp \Z) \cong \Z \oplus \varprojlim (\Z \xleftarrow{\iota^*} \Z/p_1\Z \xleftarrow{\iota^*} \Z/p_1p_2 \Z \xleftarrow{\iota^*} \Z/p_1p_2p_3 \Z \xleftarrow{\iota^*}\cdots)/ \Z,\] where each $\iota^*$ is coset inclusion.  
\end{prop}

\subsection{Finite Summability for Odometers}
Since odometers are metrically equicontinuous, we can use \cite{HSWZ} to produce unbounded cycles on $C(X) \rtimes_\vp \Z$ from the cycles on $C(X)$ constructed in the previous section. To do so, we need to extend the Hilbert Space, representation, and the operator. To start, observe that: 
\begin{prop}{(\cite[Formula 2.9]{HSWZ}, for example)}\label{hat}
    Let $A$ be a $C^*$-algebra, $H$ a separable Hilbert space, and $\pi:A \to B(H)$ a representation of $A$ on $H.$ Suppose $\vp: A \to A$ is an automorphism. Then there is a covariant representation defined on elementary tensors as follows: $\hat{\pi}: A \rtimes_\vp \Z \to B(\ell^2(\Z) \otimes H)$ via \[\hat{\pi}(a)(e_m \otimes \de_n) = e_m \otimes \pi(\vp^{-m}(a)) \de_n \text{ and } \hat{\pi}(u) (e_m \otimes \de_n) = e_{m+1} \otimes \de_n,\]where $u$ is the unitary implementing $\vp.$ 
\end{prop}  
Moving forward, if $\pi$ is a representation of $C(X),$ we denote by $\hat{\pi}$ the representation of $C(X) \rtimes_\vp \Z$ obtained by the process in Proposition \ref{hat}. In \cite{HSWZ}, in the case of metrically equicontinuous actions, the authors take a spectral triple on the base algebra and produce spectral triples on the crossed product algebra for such an action using this process to extend the representation. The results of \cite{HSWZ} that we use are summarized in: 
\begin{prop}{\cite[Section 2.4, Theorem 2.7, 2.14, and comment before Theorem 2.11]{HSWZ} }\label{hawkins}
    Let $\mathscr{A} \sub C(X)$ be dense and $\left( H_+ \oplus H_-, \rho_+ \oplus \rho_- , \begin{bmatrix} 0 & D \\ D^* & 0 \end{bmatrix} \right)$ an even spectral triple on $\A$. Let $M_\iota$ be the self-adjoint unbounded-operator in $\ell^2(Z)$ with domain $\C\Z$ and given by $M_\iota(f(n)):= nf(n).$ Then the spectral triple $(\C\Z, \ell^2(\Z), M_\iota)$ on $C^*(\Z)$ is $p$-summable for all $p >1.$ Further, if $\Z$ acts metrically equicontinuously on $X$ then:
    \begin{enumerate}
    \item $\left(H_+ \otimes \ell^2(Z) \oplus H_- \otimes \ell^2(\Z), \hat{\rho_+} \oplus \hat{\rho_-}, \begin{bmatrix} 1 \otimes M_i & D \otimes 1 \\ D^* \otimes 1 & -1 \otimes M_i \end{bmatrix} \right),$ is an odd spectral triple on $C_c(\Z, \mathscr{A}),$
    \item the spectral triple  of (i) is $q$ summable for all $q > p + 1,$
    \item and the class in $K^1(C(X) \rtimes_\vp \Z)$ of the spectral triple in (i) represents the image of $\\ \left[H_+  \oplus H_- , \rho_+ \oplus \rho_-, \begin{bmatrix} 0 & D \\ D^*  & 0\end{bmatrix} \right] \in K^0(C(X))$ under the boundary map in the Pimsner-Voiculescu exact sequence.
    \end{enumerate} 
\end{prop}
We apply this result to obtain:
\begin{theorem}\label{spectripodometer}
Let $\vp: X \to X$ be a metrically equicontinuous action on the Cantor Set (i.e. an odometer, Proposition \ref{equicontinuosodometers}) and also use $\vp:C(X) \to C(X)$ to denote the induced automorphism $f \mapsto f \circ \vp^{-1}.$  Let $\A = \text{span}\{\chi_{C_\mu}| \mu \in Y \}.$ Then, for every $x \in K^*(C(X) \rtimes_\vp \Z)$ and each $p >1,$ there is an unbounded Fredholm module representing $x$ that is $p$-summable on $C_c(\Z,\A).$
\end{theorem}
\begin{proof}
In the case that the homeomorphism $\vp$ is metrically equicontinuous, i.e. that the dynamical system is an odometer, we can use the results of \cite{HSWZ} from Proposition \ref{hawkins}. This result allows us to attain a representative for each element in $K^1(C(X) \rtimes_\vp \Z)$ as the image of the cycles from Corollary \ref{UnboundedCantor} under the Pimsner-Voiculescu boundary map. Because $\vp$ preserves $\A$ and, for each $f \in \A, \: \sup_{n \in \Z} \|[D,\pi_{\tau_\mu}(\vp^{n}(f))] \| < \infty,$ (which is guaranteed by metric equicontinuity \cite[Proposition 3.1]{HSWZ}), Proposition \ref{hawkins} gives that each class in $K^1(C(X) \rtimes_\vp \Z)$ can be represented by a cycle of the form: \[\left ( \ell^2(\Z) \otimes \ell^2(Y) \otimes \ell^2(\N) \oplus \ell^2(\Z) \otimes \ell^2(Y) \otimes \ell^2(\N), \pi_\tau', \begin{bmatrix} 1 \otimes M_i & D \otimes 1 \\ D^* \otimes 1 & -1 \otimes M_i \end{bmatrix} \right)\] where 
\begin{enumerate} 
\item $D = D^*(e_n \otimes e_\mu) = e^{|n| +|\mu|} e_n \otimes e_\mu,$ 
\item $M_i: \ell^2(\Z) \to \ell^2(\Z)$ is given by $M_i (f(n)):= n \cdot f(n), $
\item  and the representation of $C(X) \rtimes_\vp \Z$ on $\ell^2(\Z) \otimes \ell^2(Y) \otimes \ell^2(\N) \oplus \ell^2(\Z) \otimes \ell^2(Y) \otimes \ell^2(\N)$ is given by:
\begin{enumerate}
\item $\pi_\tau'(f) (e_n \otimes \xi \oplus e_m \otimes \eta) = (e_n \otimes \hat{\pi_{\tau_+}}(\vp^{-n}(f)) (\xi \oplus e_m) \otimes \hat{\pi_{\tau_-}}(\vp^{-m}(f))(\eta)$ and 
\item $ \pi_\tau'(u)(e_n \otimes \xi \oplus e_m \otimes \eta) = e_{n+1} \otimes \xi \oplus e_{m+1} \otimes \eta.$ \end{enumerate}
\end{enumerate}
Observe that $(\C \Z, \ell^2(\Z), M_i)$ is a $p$-summable spectral triple for all $p >1$, as the eigenvalues of $M_i$ are $n$ with multiplicity 1. Thus, the odd unbounded cycles on the crossed product will be $p$-summable for all $p > k+1$ when the cycles of Proposition \ref{UnboundedCantor} are $k$-summable \cite[Theorems 2.7, 2.14]{HSWZ}. Since, for each $k >0,$ we can choose these cycles on the Cantor set to be $k$-summable, for each $p>1,$ the cycles on the crossed product can be chosen to be $p$-summable. To see this, observe that \[\begin{bmatrix} 1 \otimes M_i & D \otimes 1 \\ D \otimes 1 & -1 \otimes M_i \end{bmatrix}^2 = \begin{bmatrix} D^2 \otimes 1 + 1 \otimes M_i^2 
 & 0 \\ 0 & D^2 \otimes 1 + 1 \otimes M_i^2 \end{bmatrix},\] so that the eigenvalues of the square are $W^{2(n+|\mu|)} + m^2$ with associated eigenvectors \[\left \{\begin{bmatrix} e_n \otimes e_\mu \otimes e_m \\ 0\end{bmatrix}, \begin{bmatrix}0 \\ e_n \otimes e_\mu \otimes e_m \end{bmatrix}\right \}.\] The eigenvalues of $\left(1+ \begin{bmatrix} 1 \otimes M_i & D \otimes 1 \\ D \otimes 1 & -1 \otimes M_i \end{bmatrix}^2 \right)^{-\f p2}$ are $(1+W^{2(n+|\mu|)} + m^2)^{-\f p2}$ with the same eigenvectors, so that \[\sum_{n \in \Z,\: \mu \in Y, \: m \in \Z} (1+W^{2(n+|\mu|)} + m^2)^{-\f p2} \] converges if and only if $p \geq 1  + k$ where $D$ is $k$-summable.
\end{proof}
\begin{Cor}
    Under the bounded transform $(H, \rho, D) \mapsto \left(H, \rho, F =\frac{D}{(1+D^2)^{\f 12}}\right),$ the unbounded cycles of the previous proposition \[\left (\ell^2(\Z) \otimes \ell^2(Y) \otimes \ell^2(\N) \oplus \ell^2(\Z) \otimes \ell^2(Y) \otimes \ell^2(\N), \pi_\tau', \begin{bmatrix} 1 \otimes M_i & D \otimes 1 \\ D^* \otimes 1 & -1 \otimes M_i \end{bmatrix} \right)\] map to bounded cycles that are $p$-summable  for all $p >1$ on $C_c(\Z,\A),$ where $\A$ is defined in Theorem \ref{spectripodometer}.
\end{Cor}
\begin{proof}
Observe that, \[F = \begin{bmatrix} 1 \otimes M_i & D \otimes 1 \\ D \otimes 1 & -1 \otimes M_i \end{bmatrix} \left (1 + \begin{bmatrix} 1 \otimes M_i & D \otimes 1 \\ D \otimes 1 & -1 \otimes M_i \end{bmatrix}^2 \right)^{-\f 12} = (1 + D^2 \otimes 1 + 1 \otimes M_i^2)^{-\f 12} \begin{bmatrix} 1 \otimes M_i & D \otimes 1 \\ D \otimes 1 & -1 \otimes M_i \end{bmatrix},\] so that \[\begin{bmatrix} 1 \otimes M_i & D \otimes 1 \\ D \otimes 1 & -1 \otimes M_i \end{bmatrix} \left (1 + \begin{bmatrix} 1 \otimes M_i & D \otimes 1 \\ D \otimes 1 & -1 \otimes M_i \end{bmatrix}^2 \right)^{-\f 12} \begin{bmatrix}e_n \otimes e_\mu \otimes e_m \\ e_n \otimes e_\mu \otimes e_m \end{bmatrix} = \begin{bmatrix} \frac{m + W^{n+|\mu|}}{\sqrt{1 + W^{2(n+|\mu|)}}} e_n \otimes e_\mu \otimes e_m \\ \frac{W^{n+|\mu|}-m}{\sqrt{1 + W^{2(n+|\mu|)}}} e_n \otimes e_\mu \otimes e_m \end{bmatrix}.\] Now, evaluating the commutators \[[F,\pi_\tau'(\chi_{C_\mu})] \left (\begin{bmatrix} e_n \otimes e_\mu \otimes e_m \\ e_n \otimes e_\mu \otimes e_m\end{bmatrix} \right) = \] \[\begin{bmatrix} \frac{m}{\sqrt{1 +m^2 W^{2(n+|\mu|)}}} e_n \otimes \pi_{\tau_{n_+}} (\vp^{-m}((\chi_{C_\mu})) e_\mu \otimes e_m + \frac{W^{n+|\mu|)}}{\sqrt{1+m^2 + W^{2(n+|\mu|)}}}  e_n \otimes (\pi_{\tau_{n_-}}(\vp^{-m}(\chi_{C_\mu})) e_\mu \otimes e_m \\  \frac{W^{n+|\mu|}}{\sqrt{1 +m^2 + e^{2(n+|\mu|)}}}  e_n \otimes \pi_{\tau_{n_+}} (\vp^{-m}((\chi_{C_\mu})) - \frac{m}{\sqrt{1 +m^2+ W^{2(n+|\mu|)}}}  e_n \otimes (\pi_{\tau_{n_-}}(\vp^{-m}(\chi_{C_\mu})) e_\mu \otimes e_m\end{bmatrix}\] \[- \begin{bmatrix} \pi_{\tau_{n_+}} (\vp^{-m}(\chi_{C_\mu})) \frac{m + W^{n+|\mu|}}{\sqrt{1 +m^2 + W^{2(n+|\mu|)}}}e_n \otimes e_\mu \otimes e_m \\ \pi_{\tau_{n_-}}(\vp^{-m}(\chi_{C_\mu})) \frac{W^{n+|\mu|}-m}{\sqrt{1 + m^2+ W^{2(n+|\mu|}}}   e_n \otimes e_\mu \otimes e_m\end{bmatrix} \] \[ = \begin{bmatrix} \frac{W^{n+|\mu|}}{\sqrt{1+m^2+W^{2(n+|\mu|)}}} (\pi_{\tau_{n_+}}(\vp^{-m}(\chi_{C_\mu}))-(\pi_{\tau_{n_-}}(\vp^{-m}(\chi_{C_\mu})) e_n \otimes e_\mu \otimes e_m \\ \frac{-W^{n+|\mu|}}{\sqrt{1+m^2+W^{2(n+|\mu|)}}} (\pi_{\tau_{n_+}}(\vp^{-m}(\chi_{C_\mu}))-(\pi_{\tau_{n_-}}(\vp^{-m}(\chi_{C_\mu})) e_n \otimes e_\mu \otimes e_m  \end{bmatrix}.\]
Note that, for each $m$, 
\[
\sum_{n \in \N} \left( \pi_{\tau_{n_+}}(\vp^{-m}(\chi_{C_\mu})) - \pi_{\tau_{n_-}}(\vp^{-m}(\chi_{C_\mu})) \right) 
= \gm(x)([\chi_{C_{\vp^m(C_\mu)}}]),
\]
for $x$ the class represented by the initial spectral triple on $C(X).$ Because the odometer permutes cylinder sets of the same length, the orbit of $C_\mu$ returns to itself. Thus, the quantity 
$k := \max\{\gm(x)([\chi_{\vp^m(C_\mu)}]) \mid m \in \Z\}$
is finite. Further, $\pi_{\tau_{n_+}}(\vp^{-m}((\chi_{C_\mu}))$ and $\pi_{\tau_{n_-}}(\vp^{-m}((\chi_{C_\mu}))$ only differ on words of length less than $|\mu|$ and for finitely many $n.$ Thus, there is an $M < \infty$, so that
\[
\left\|[F,\pi_\tau'(\chi_{C_\mu})]
\begin{bmatrix} 
e_n \otimes e_\mu \otimes e_m \\ 
e_n \otimes e_\mu \otimes e_m
\end{bmatrix}
\right\| 
\leq \frac{M}{\sqrt{1+m^2}}
\begin{bmatrix} 
(\pi_{\tau_{n_+}}(\vp^{-m}(\chi_{C_\mu})) - \pi_{\tau_{n_-}}(\vp^{-m}(\chi_{C_\mu}))) e_n \otimes e_\mu \otimes e_m \\ 
(\pi_{\tau_{n_+}}(\vp^{-m}(\chi_{C_\mu})) - \pi_{\tau_{n_-}}(\vp^{-m}(\chi_{C_\mu}))) e_n \otimes e_\mu \otimes e_m 
\end{bmatrix}.
\]
\end{proof}
\bibliographystyle{plain}	
\bibliography{refs}

\end{document}